\documentclass[11pt]{amsart}

\usepackage{amsopn, amssymb, amscd, amsmath, amsthm, MnSymbol, stmaryrd}
\usepackage{enumerate}
\usepackage{color}
\usepackage [all] {xy}
\usepackage[perpage,symbol*]{footmisc}
\usepackage{array}
\usepackage{multirow}
\usepackage{graphicx}
\usepackage{color, colortbl}
\usepackage{array, booktabs, arydshln, xcolor}
\usepackage{dsfont}
\usepackage{mathrsfs}
\usepackage{anysize}
\marginsize{3cm}{3cm}{3cm}{3cm}

\newcommand{\nc}{\newcommand}

\nc{\gramU}{\textbf{\textsl{U}}}

\nc{\fg}{\mathfrak{f} }     \nc{\vg}{\mathfrak{v} }       \nc{\wg}{\mathfrak{w} }
\nc{\zg}{\mathfrak{z} }     \nc{\ngo}{\mathfrak{n} }      \nc{\kg}{\mathfrak{k} }
\nc{\ngoc}{\widehat{\mathfrak{n}} }
\nc{\mg}{\mathfrak{m} }     \nc{\bg}{\mathfrak{b} }       \nc{\ggo}{\mathfrak{g} }
\nc{\ggoc}{\widehat{\mathfrak{g}} }
\nc{\sog}{\mathfrak{so} }
\nc{\sug}{\mathfrak{su} }   \nc{\spg}{\mathfrak{sp} }     \nc{\slg}{\mathfrak{sl} }
\nc{\glg}{\mathfrak{gl} }   \nc{\cg}{\mathfrak{c} }       \nc{\rg}{\mathfrak{r} }
\nc{\hg}{\mathfrak{h} }     \nc{\tg}{\mathfrak{t} }       \nc{\ug}{\mathfrak{u} }
\nc{\dg}{\mathfrak{d} }     \nc{\ag}{\mathfrak{a} }       \nc{\pg}{\mathfrak{p} }
\nc{\agc}{\widehat{\mathfrak{a}} }
\nc{\sg}{\mathfrak{s} }     \nc{\affg}{\mathfrak{aff} }
\nc{\ggob}{\overline{\mathfrak{g}} }

\nc{\pca}{\mathcal{P}}       \nc{\nca}{\mathcal{N}}       \nc{\lca}{\mathcal{L}}
\nc{\oca}{\mathcal{O}}       \nc{\mca}{\mathcal{M}}       \nc{\tca}{\mathcal{T}}
\nc{\aca}{\mathcal{A}}       \nc{\cca}{\mathcal{C}}       \nc{\gca}{\mathcal{G}}
\nc{\sca}{\mathcal{S}}       \nc{\hca}{\mathcal{H}}       \nc{\bca}{\mathcal{B}}
\nc{\dca}{\mathcal{D}}       \nc{\rca}{\mathcal{R}}

\nc{\val}{\operatorname{val}}

\nc{\vp}{\varphi}
\nc{\ddt}{\tfrac{{\rm d}}{{\rm d}t}}
\nc{\dpar}{\tfrac{\partial}{\partial t}}

\nc{\im}{\sqrt{-1}}        

\newcommand{\comillas}[1]{\textquotedblleft{#1}\textquotedblright}

\nc{\SO}{\mathrm{SO}}           \nc{\Spe}{\mathrm{Sp}}          \nc{\Sl}{\mathrm{SL}}
\nc{\SU}{\mathrm{SU}}           \nc{\Or}{\mathrm{O}}            \nc{\U}{\mathrm{U}}
\nc{\Se}{\mathrm{S}}            \nc{\Cl}{\mathrm{Cl}}           \nc{\Spein}{\mathrm{Spin}}
\nc{\Pin}{\mathrm{Pin}}
\nc{\Glr}{\mathrm{GL}_n(\RR)}   \nc{\Glc}{\mathrm{GL}_n(\CC)}   \nc{\Glv}{\mathrm{GL}(V)}    \nc{\Glk}{\mathrm{GL}_n(\fk)}   \nc{\Gl}{\mathrm{GL}}
\nc{\GrpG}{\mathrm{G}}          \nc{\GrpH}{\mathrm{H}}          \nc{\GrpA}{\mathrm{A}}       \nc{\GrpT}{\mathrm{T}}          \nc{\GrpK}{\mathrm{K}}
\nc{\GrpGc}{\widehat{\mathrm{G}}}

\nc{\g}{\mathfrak{gl}_n(\RR)}

\nc{\RR}{{\Bbb R}} \nc{\HH}{{\Bbb H}} \nc{\CC}{{\Bbb C}} \nc{\ZZ}{{\Bbb Z}}
\nc{\FF}{{\Bbb F}} \nc{\NN}{{\Bbb N}} \nc{\QQ}{{\Bbb Q}} \nc{\PP}{{\Bbb P}}

\nc{\euler}{{\rm e}}

\nc{\vs}{\vspace{.2cm}} \nc{\vsp}{\vspace{1cm}}

\nc{\ip}{\langle \cdot , \cdot \rangle}
\nc{\ipd}{\langle \hspace{-0.5mm}\langle \cdot , \cdot \rangle\hspace{-0.5mm}\rangle}
\nc{\ipp}{(      \cdot , \cdot       )}
\nc{\la}{\langle} \nc{\ra}{\rangle}

\nc{\ortsum}{ \mbox{\tiny $\displaystyle \bigoplus^{\perp}$}}

\nc{\dirsum}{ \mbox{\tiny $\displaystyle \bigoplus $}}

\nc{\ortres}{ \mbox{\tiny $\displaystyle \bigominus^{\perp}$}}

\nc{\unm}{\tfrac{1}{2}}\nc{\unc}{\tfrac{1}{4}} \nc{\und}{\tfrac{1}{16}}

\nc{\no}{\vs\noindent}

\nc{\lamn}{\Lambda^2(\RR^n)^*\otimes\RR^n} \nc{\lamp}{\Lambda^2\pg^*\otimes\pg}
\nc{\lamg}{\Lambda^2\ggo^*\otimes\ggo} \nc{\lamngo}{\Lambda^2\ngo^*\otimes\ngo}
\nc{\lamnk}{\Lambda^2(\fk^n)^*\otimes \fk^n} \nc{\lamnkt}{\Lambda^3(\fk^n)^*\otimes \fk^n}

\nc{\tangz}{{\rm T}^{\rm Zar}}
\nc{\mum}{/\!\!/} \nc{\kir}{/\!\!/\!\!/}
\nc{\lievark}{\mathfrak{L}_n(\fk)}         \nc{\lievarc}{\mathfrak{L}_n(\CC)}        \nc{\lievarr}{\mathfrak{L}_n(\RR)}
\nc{\solvvark}{\mathfrak{R}_n(\fk)}        \nc{\solvvarc}{\mathfrak{R}_n(\CC)}       \nc{\solvvarr}{\mathfrak{R}_n(\RR)}
\nc{\nilvark}{\mathfrak{N}_n(\fk)}         \nc{\nilvarc}{\mathfrak{N}_n(\CC)}        \nc{\nilvarr}{\mathfrak{N}_n(\RR)}
\nc{\cirre}{\textrm{C}}
\nc{\fk}{\mathrm{k}}

\nc{\Ri}{\tfrac{4\Ric_{\mu}}{||\mu||^2}}

\nc{\ds}{\displaystyle}

\nc{\lb}{[\cdot,\cdot]}

\nc{\Hess}{\operatorname{Hess}}
\nc{\diag}{\operatorname{Diag}}   \nc{\Id}{\operatorname{Id}}
\nc{\trans}{\mbox{{\tiny$\operatorname{T}$}}}                         \nc{\Proj}{\operatorname{Proj}}
\nc{\Proy}{\operatorname{Proy}}
\nc{\ad}{\operatorname{ad}}       \nc{\Ad}{\operatorname{Ad}}        
\nc{\rank}{\operatorname{rank}}   \nc{\codim}{\operatorname{codim}}  
\nc{\Irr}{\operatorname{Irr}}     \nc{\End}{\operatorname{End}}
\nc{\Aut}{\operatorname{Aut}}     \nc{\Inn}{\operatorname{Inn}}
\nc{\lRad}{\operatorname{Rad}}
\nc{\Der}{\operatorname{Der}}     \nc{\Ker}{\operatorname{Ker}}
\nc{\Iso}{\operatorname{I}}       \nc{\Diff}{\operatorname{Diff}}
\nc{\Lie}{\operatorname{Lie}}     \nc{\tr}{\operatorname{tr}}
\nc{\dif}{\operatorname{d}}       \nc{\e}{\operatorname{e}}
\nc{\sen}{\operatorname{sen}}     \nc{\tang}{\operatorname{T}}
\nc{\modu}{\operatorname{mod}}
\nc{\Riem}{\operatorname{R}}     \nc{\Ric}{ {Ric}}
\nc{\curvature}{\operatorname{K}}
\nc{\sym}{\operatorname{sym}}     \nc{\symac}{\operatorname{sym^{ac}}}   \nc{\symc}{\operatorname{sym^{c}}}
\nc{\scalar}{\operatorname{sc}}
\nc{\grad}{\operatorname{grad}}
\nc{\ricci}{\operatorname{ric}}   \nc{\nr}{\operatorname{nr}}            \nc{\riccic}{\operatorname{ric^{c}}}
\nc{\riccig}{\operatorname{ric^{\gamma}}}
\nc{\Rin}{\operatorname{M}}
\nc{\Kill}{\operatorname{B}}
\nc{\Le}{\operatorname{L}}
\nc{\level}{\operatorname{level}} \nc{\rad}{\operatorname{r}}
\nc{\abel}{\operatorname{ab}}
\nc{\CH}{\operatorname{CH}}        \nc{\mcc}{\operatorname{mcc}}     \nc{\inte}{\operatorname{int}}         \nc{\aff}{\operatorname{Aff}}
\nc{\CaC}{\operatorname{CC}}        \nc{\ccm}{\operatorname{ccm}}
\nc{\Adj}{\operatorname{Adj}}
\nc{\Order}{\operatorname{O}} \nc{\Ricg}{\operatorname{Ric^{\gamma}}}
\nc{\Hom}{\operatorname{Hom}}
\nc{\sign}{\operatorname{sign}}
\nc{\spanv}{\operatorname{span}}
\nc{\xp}{\operatorname{xp}}   \nc{\xt}{\operatorname{xt}}
\nc{\IC}{\operatorname{IC}}   \nc{\OC}{\operatorname{OC}}

\nc{\rhov}{\operatorname{\rho_{v}}}
\nc{\mm}{m}
\nc{\mmt}{\widetilde{m}}
\nc{\F}{\operatorname{F}}

\theoremstyle{plain}
\newtheorem{theorem}{Theorem}[section]

\theoremstyle{definition}

\newtheorem{notation}[theorem]{Notation}

\theoremstyle{remark}
\newtheorem{remark}[theorem]{Remark}

\newtheorem{example}[theorem]{Example}

\numberwithin{equation}{section}

\begin{document}

\title{Classification of Nilsoliton metrics in dimension seven.}

\author{EDISON ALBERTO FERN\'ANDEZ CULMA}

\address{Current affiliation: CIEM, FaMAF, Universidad Nacional de C\'ordoba, \newline \indent Ciudad Universitaria, \newline \indent (5000) C\'ordoba, \newline \indent Argentina}

\email{efernandez@famaf.unc.edu.ar}

\thanks{Fully supported by a CONICET postdoctoral fellowship (Argentina)}



\maketitle
\begin{abstract}
The aim of this paper is to classify Ricci soliton metrics on $7$-dimensional nilpotent Lie groups.
It can be considered as a continuation of our paper in Transformation Groups, Volume 17, Number 3 (2012), 639--656.
To this end, we use the classification of $7$-dimensional real nilpotent Lie algebras given
by Ming-Peng Gong in his Ph.D thesis and some techniques from the results of Michael Jablonski in \cite{JABLONSKI2,JABLONSKI1} and
of Yuri Nikolayevsky in \cite{NIKOLAYEVSKY2}.
Of the $9$ one-parameter families and $140$ isolated $7$-dimensional indecomposable real nilpotent Lie algebras, we have $99$ nilsoliton metrics given
in an explicit form and $7$ one-parameter families admitting nilsoliton metrics.

Our classification is the result of a case-by-case analysis, so many illustrative examples are carefully developed to explain how to obtain the main result.

\end{abstract}

\section{Introduction}
The Ricci flow, introduced by Richard Hamilton in the early 1980s,  is a geometric evolution equation which deforms smoothly an initial metric $g_0$ on a Riemannian manifold $M$ in the direction of minus two times its Ricci tensor:
$$\left\{\begin{array}{l}
  \dpar g = -2 \ricci_{g},\\ g(0)=g_{0}.
  \end{array}\right.
$$
The \comillas{philosophy} behind this tool of geometric analysis is to try to evolve the geometry of $(M,g_0)$
to one which looks \comillas{more uniform} (although in practice, it is not always possible).

A \textit{Ricci soliton} is a complete Riemannian metric $g_{0}$ such that the solution to the Ricci flow $g(t)$ with $g(0)=g_{0}$ changes $g_{0}$ only by  diffeomorphisms and scaling as time goes on, that is
$g(t)=c(t)\varphi_{t}^{\ast}g_{0}$, where $c(t)\in\RR_{+}$ and $\varphi_{t}$ is a one-parameter group
of diffeomorphisms of $M$; a Ricci soliton is not \comillas{improved} by the Ricci flow.
These \comillas{distinguished} metrics are important in the study of the Ricci flow because they may be
limiting cases for the Ricci flow near \textit{singularities} and are a natural generalization of an \textit{Einstein metric} (\cite{CAO}).

{

In Ricci Flow theory, homogeneous Riemannian manifolds provide a rich source of explicit examples for some concepts and
behaviors. Moreover, many of the known examples of Ricci solitons correspond to \textit{algebraic Ricci solitons}: these are given by a left-invariant metric $g$ on a simply-connected Lie group $\GrpG$ satisfying
\begin{eqnarray}
\Ric(g)=c\Id+D,\, \mbox{for some } c\in\RR \mbox{ and } D\in \Der(\ggo)
\end{eqnarray}
where $\ggo=\Lie(\GrpG)$ and $\Ric$ is the Ricci operator of $g$.

Such metrics are indeed Ricci solitons and in the particular case that $\GrpG$ is solvable (respectively nilpotent) are called \textit{solvsolitons} (respectively \textit{nilsolitons}).

Our aim in this paper is to classify nilsolitons on $7$-dimensional simply connected nilpotent Lie groups.
In general, it is difficult to know if a nilpotent Lie algebra admits a nilsoliton (inner product) and
it is very difficult to get such metric when it exists.
Nilsoliton metrics have been completely classified only up to dimension $6$ by Jorge Lauret and Cynthia Will in \cite{LAURET2, WILL1} (with a minor correction in \cite{LAURETW2}). In \cite{PAYNE1}, Trayce Payne and H\"{u}lya Kadio\u{g}lu have introduced a computational method for classifying nilsoliton metrics in the family of nilpotent Lie algebras with \textit{simple pre-Einstein derivation} and \textit{nonsingular Gram matrix}. Their method does not rely on any preexisting classifications of nilpotent Lie algebras. They classify nilsolitons within this family, which has $33$ algebras in dimensions $7$ and $159$ algebras in dimension $8$. In \cite{FERNANDEZ-CULMA1}, we gave a complete classification of all seven-dimensional (complex)
nilpotent Lie algebras $\ngo$ such that a real form of $\ngo$ admits a nilsoliton inner product.  In order, to prove  \cite[Theorem 7]{FERNANDEZ-CULMA1}, we make heavy use of results proved by Yuri Nikolayevsky in \cite{NIKOLAYEVSKY2} (some of such results were independently proved by Michael Jablonski in \cite{JABLONSKI1, JABLONSKI2} in a more general form).

In these notes, we use the classification of seven-dimensional real nilpotent Lie algebras given by Ming-Peng Gong in his Ph.D thesis (\cite{GONG1}), results from
complex (and real) Geometric Invariant Theory (\cite{NESS1,JABLONSKI1, JABLONSKI2,FERNANDEZ-CULMA3}) and calculations given in \cite{FERNANDEZ-CULMA2}.

Following \cite{GONG1}, the indecomposable $7$-dimensional real Lie algebras can be seen as $9$ one-parameter families of nilpotent Lie algebras plus $140$ \comillas{isolated} nilpotent Lie algebras; i.e. $140$ algebras that are not in any of such one-parameter families. As we mentioned before, to give a nilsoliton metric can be a very complicated problem and, in the case of one-parameter families admitting nilsoliton metrics, to give a nilsoliton metric explicitly for each member of the family in a simultaneous manner it seems almost impossible. Therefore, we focus our attention in isolated algebras and in the family $(147E)[0<t<1]$ where we can give a nilsoliton metric for each $t$. So, we have $99$ nilsoliton metrics and $7$ one-parameter families admitting nilsoliton metrics. The methods explained here can also be used to study any fixed member in a one-parameter family.

Our classification is the result of a case-by-case analysis, so many illustrative examples are carefully developed to explain how to obtain the main result.
Some nilsoliton metrics here were obtained in \cite{FERNANDEZ-CULMA2}, we refer the reader to \cite{FERNANDEZ-CULMA1} for more details.

\section{Preliminaries}

In this section, we give a brief exposition on nilsoliton metrics and real (and complex) geometric invariant theory.
There is an intriguing interplay between the Ricci flow on nilpotent Lie groups and the gradient flow of the norm squared of the \textit{moment map} associated to the natural action of $\Glr$ on $V=\lamn$. It is known that if $\ngo:=(\RR^{n},\mu)$ is a nilpotent Lie algebra and $\ip$ is the canonical inner product of $\RR^n$, then the Ricci operator of $(\RR^n,\mu,\ip)$ satisfies
\begin{eqnarray}
  4\Ric &=& \mm(\mu)
\end{eqnarray}
where $\mm$ is the moment map. It follows that to minimize the norm of the Ricci tensor among all left-invariant metrics of $\ngo$ with the same scalar curvature is equivalent to minimize $||\mm||^2/||\mu||^4$ along the $\Glr$-orbit of $\mu$ (here, the inner products on $\g$ and $V$, which are denoted by $\ipd$ and $\ip$ respectively, are those induced by the canonical inner product of $\RR^n$).

\begin{theorem}[Jorge Lauret](see for instance \cite[Theorem 4.2]{LAURET3})
The nilpotent Lie algebra $\ngo:=(\RR^n,\mu)$ admits a nilsoliton (inner product) if and only if the $\Glr$-orbit of $\mu$ is distinguished; i.e.
$\Glr \cdot \mu$ contains a critical point of $||\mm||^2$, or equivalently, if and only if there exists a $g \in \Glr$ such that
\begin{eqnarray}\label{defnil}
  \mm(\widetilde{\mu}) = c\Id + D \, \mbox{for some } c\in\RR \mbox{ and } D\in \Der(\widetilde{\ngo})
\end{eqnarray}
where $\widetilde{\mu}:=g\cdot\mu$ and $\widetilde{\ngo}:=(\RR^n,\widetilde{\mu})$.

There is at most one nilsoliton metric on a nilpotent Lie group (up to isometry and scaling).
\end{theorem}

A point $\widetilde{\mu}$ satisfying the Equation (\ref{defnil}) is called a \textit{distinguished point} (see \cite[Definition 2.6]{JABLONSKI2}) and
the derivation $D$ is usually called \textit{Einstein derivation} (in connection with \textit{Einstein solvmanifolds}) or \textit{nilsoliton derivation}.

It is know from \cite{HEBER1} that the eigenvalues of an Einstein derivation are all positive integers without a common divisor (up to a rational factor) and that it is, up to conjugation by an automorphism,  positively proportional to the \textit{pre-Einstein derivation}, which is defined by Yuri Nikolayevsky in \cite[Definiton 2]{NIKOLAYEVSKY2}.

Before stating the techniques to classify the nilsoliton metrics, let us first recall some notations.

\begin{notation}
  Let $\ag$ denote the maximal abelian subalgebra of $\g$ contained in the vector space of symmetric matrices $\pg$ given by
\begin{eqnarray*}
  \ag &=& \{diag(x_{1},\ldots,x_{n}) : x_{i} \in \RR\}
\end{eqnarray*}
and set $\GrpA:=\exp(\ag)$.

Let $\Phi$ be a finite subset of $\ag$. The convex hull of $\Phi$ will be denoted by $\CH (\Phi)$ while by $\aff(\Phi)$ we denote the affine space generated by $\Phi$. $\mcc (\Phi)$ denotes the \textit{minimal convex combination} of $\Phi$; i.e. the unique vector closest to the origin in $\CH (\Phi)$. The notation $\inte(\CH(\Phi))$ represents the interior of $\CH(\Phi)$ relative to the usual topology of $\aff(\Phi)$.

Given a nilpotent Lie algebra $\ngo:=(\RR^n,\mu)$ we denote by $\mathfrak{R}(\mu)$ the (ordered) set of weights related with $\mu$ to the action of $\Glr$ on $V$, i.e. if $\{C_{i,j}^{k}\}$ are the structural constants of $\ngo$ in the basis $\{e_1\ldots e_n\}$ then
\begin{eqnarray*}
  \mathfrak{R}(\mu) &=& \{E_{k,k}-E_{i,i}-E_{j,j} : C_{i,j}^{k}\neq 0\}
\end{eqnarray*}
where $\{E_{i,j}\}$ is the canonical basis of $\g$.
\end{notation}

Recall that Nikolayevsky's nice basis criterium (\cite[Theorem 3]{NIKOLAYEVSKY2}) says that a nilpotent Lie algebra $\ngo:=(\RR^n,\mu)$ which is written in a \textit{nice basis} (\cite[Definition 3]{NIKOLAYEVSKY2}) admits a nilsoliton metric if and only if the Equation
\begin{eqnarray}\label{eqn1}
\gramU x = [1]_{m}
\end{eqnarray}
has at least one solution $x$ with positive coordinates, where $m=\#(\mathfrak{R}(\mu))$ and $\gramU$ is the Gram matrix of $(\mathfrak{R}(\mu),\ipd)$. This result give us an easy-to-check convex geometry condition for a nilpotent Lie algebra with a nice basis to be an Einstein nilradical.

If one carefully reads the proof of the above theorem one sees that it gives us a technique to find a nilsoliton on a nilpotent Lie algebra which is written in a nice basis and admits such distinguished metric; a nilsoliton can be found it in a $\GrpA$-orbit. So, Nikolayevsky's nice basis criterium can be rewritten as

\begin{theorem}\label{niceTh}[Nikolayevsky's nice basis criterium]
Let $\ngo:=(\RR^n,\mu)$ be a nilpotent Lie algebra such that $\ngo$ is written in a nice basis. $\ngo$ admits a nilsoliton if and only if
there exist $g\in\GrpA$ such that
\begin{eqnarray}\label{eqn2}
  \mm(g\cdot\mu)&=&\mcc(\mathfrak{R}(\mu)) ,
\end{eqnarray}
and, consequently, $\widetilde{\mu}:=g\cdot\mu$ is a distinguished point in the $\Glr$-orbit of $\mu$.
\end{theorem}

To find a nilsoliton metric for a nilpotent Lie algebra admitting such metric and which is written in a nice basis is easy in practice. We must calculate the vector $\mcc(\mathfrak{R}(\mu))$, which is given by
$$
\frac{1}{\sum x_p}\left( \sum x_p \mathfrak{R}_{}(\mu)_{p} \right)
$$
where $[x_i]$ is any positive solution to the Equation (\ref{eqn1}) and we solve the Equation (\ref{eqn2}) for $g\in \GrpA$. We refer the reader to \cite[Corollary 3.4]{JABLONSKI1} or \cite[Section 3]{FERNANDEZ-CULMA3} for further information on results related with Theorem \ref{niceTh}.

Another application from geometric invariant theory in the study of nilsoliton metrics is the following result, which was proven independently Michael Jablonski in \cite[Theorem 6.5]{JABLONSKI1} and Yuri Nikolayevsky in \cite[Theorem 6]{NIKOLAYEVSKY2}.

\begin{theorem}\label{complex}
Let $\ngo_{1}$ and $\ngo_{2}$ be two (real) nilpotent Lie algebras whose complexifications are isomorphic as complex nilpotent Lie algebras. If $\ngo_{1}$ is an Einstein nilradical then so is $\ngo_{2}$, with the same eigenvalue type.
\end{theorem}

Consider the natural action of the complex reductive Lie group $\mathrm{GL}_n(\CC)$ on $\Lambda^2(\CC^n)^*\otimes\CC^n$ and its moment map $\check{\mm}$ as in \cite{NESS1} (see \cite{LAURET4}). The Theorem \ref{complex} can be derived of comparing the distinguished orbits of the mentioned action with the distinguished orbits of the natural action of $\mathrm{GL}_n(\RR)$ on $\Lambda^2(\RR^n)^*\otimes\RR^n$. We can rephrase the Theorem \ref{complex} as saying that

\begin{theorem}\label{complex1}\cite[Theorem 4.7]{JABLONSKI1}
Let $\ngo:=(\RR^n,\mu)$ be a nilpotent Lie algebra. $\ngo$ admits a nilsoliton metric if and only if the $\mathrm{GL}_n(\CC)$-orbit of $\mu$ is distinguished of the action of $\mathrm{GL}_n(\CC)$ on $\Lambda^2(\CC^n)^*\otimes\CC^n$.
\end{theorem}

 Theorem \ref{complex1} provides us with another technique. It is easy to see that Nikolayevsky's nice basis criterium is also true in the complex case (see \cite[Remark 3.2]{FERNANDEZ-CULMA3}); i.e. given a complex nilpotent Lie algebra $(\CC^n,\mu)$ which is written in a nice basis for the action of $\mathrm{GL}_n(\CC)$ on $\Lambda^2(\CC^n)^*\otimes\CC^n$  ($\check{\mm}(\GrpA \cdot \mu) \subseteq \ag$), then $\mathrm{GL}_n(\CC) \cdot \mu$ is distinguished if and only if the Equation \ref{eqn1} has at least one solution $x$ with positive coordinates where $m=\#(\mathfrak{R}(\mu))$ and $\gramU$ is the Gram matrix of $(\mathfrak{R}(\mu),\ipd)$ (here, $\ipd$ is the usual Hermitian inner product on $\mathfrak{gl}_n(\CC)$).

A real nilpotent Lie algebra $\ngo:=(\RR^n,\mu)$ could not admit a nice basis for the action of $\mathrm{GL}_n(\RR)$ on $\Lambda^2(\RR^n)^*\otimes\RR^n$. However, it may happen that $\ngo$ admits a nice basis for the action of $\mathrm{GL}_n(\CC)$ on $\Lambda^2(\CC^n)^*\otimes\CC^n$. Suppose that $(\CC^n,\widehat{\mu})$ is written in a nice basis with $\widehat{\mu} \in \mathrm{GL}_n(\CC)\cdot\mu$ and $\mathrm{GL}_n(\CC)\cdot\mu$ being a distinguished orbit; consequently, $(\RR^n,\mu)$ admitting a nilsoliton metric. One can easily find a distinguished point $\widetilde{\mu}$ in the $\GrpA$-orbit of $\widehat{\mu}$ as above. By results of Linda Ness \cite[Theorem 6.2]{NESS1}, it is well known that any distinguished point in $\mathrm{GL}_n(\CC)\cdot\mu$ is in the $\CC^{\ast}\U(n)$-orbit of $\widetilde{\mu}$, hence to find a distinguished point in the $\mathrm{GL}_n(\RR)$-orbit of $\mu$, we study the real forms in $\CC^{\ast}\U(n)\cdot\widetilde{\mu}$ which are isomorphic over $\RR$ to $\mu$.

\section{The Classification}
In this section, we give the classification of nilsoliton metrics on $7$-dimensional (isolated) nilpotent Lie algebras. By using that any nilpotent Lie algebra of dimension less or equal than $6$ is an Einstein nilradical, one obtains that any decomposable $7$-dimensional nilpotent Lie algebra is an Einstein nilradical and it is easy to give a nilsoliton metric in each case.

Following to \cite{GONG1}, there are $5$-one parameter families and $41$ isolated nilpotent Lie algebras which are not written in a nice basis.
We will use the exclamation mark to indicate such algebras, and where we have an extra exclamation mark, such algebra has a non-positive pre-Einstein derivation; and hence it does not admit any nilsoliton metric. The one-parameter families $(147E1)[t]$, $(1357S)[t]$ and the isolated algebras $(257J1)$, $(247E)$, $(247G)$, $(247H)$, $(247H1)$, $(247R)$, $(1357Q)$, $(1357Q1)$, $(1357R)$, $(12457L)$ can be worked out as the Examples \ref{ex2} and \ref{ex3},  while the remaining algebras were studied in \cite{FERNANDEZ-CULMA2} (to obtain \cite[Theorem 7]{FERNANDEZ-CULMA1}, we give in some cases nilsoliton metrics which are considered here again).

\begin{theorem}
The classification of $7$-dimensional nilsoliton metrics on isolated nilpotent Lie algebras (plus $(147E)[0<t<1]$) is given according to the notation in \cite{GONG1} by the following list
\end{theorem}

The notation in the list is as follows: $\dim(\Der):=\mbox{Dimension of the algebra of derivations}$, $\rank:=\mbox{Dimension of a maximal torus of semisimple derivations}$. By $\phi$ we denote $4$ times the Einstein derivation or the pre-Einstein derivation, depending if the respective algebra admits or does not admit a nilsoliton metric, and $||\beta||^2$ is the real number such that
\begin{eqnarray*}
  4\Ric &=& -||\beta||^2\Id + \phi
\end{eqnarray*}
Each Lie bracket $\widetilde{\mu}$ given in the list is such that the canonical inner
product of $\RR^n$ is a nilsoliton metric of scalar curvature equal to $-\frac{1}{4}$
on $\widetilde{\ngo}:=(\RR^7,\widetilde{\mu})$.

\begin{enumerate}

\item \textbf{(37A)}: $\dim(\Der)=25$, $\rank=4$, Betti Numbers $(4, 12, 18, 18, 12, 4, 1)$ ($\cong_{\RR} \ggo_{4.2}$).
  \newline Einstein der. $\phi=\frac{2}{3}\diag(2,1,2,2,3,3,3)$, $||\beta||^2=\frac{5}{3}\thickapprox1.667$
  \newline $ [{e_1},{e_2}]=\frac{\sqrt{6}}{6}{e_5},[{e_2},{e_3}]=\frac{\sqrt{6}}{6}{e_6},[{e_2},{e_4}]=\frac{\sqrt{6}}{6}{e_7} $

\item \textbf{(37B)}: $\dim(\Der)=20$, $\rank=4$, Betti Numbers $(4, 11, 16, 16, 11, 4, 1)$ ($\cong_{\RR} \ggo_{4.1}$).
  \newline Einstein der. $\phi=\frac{1}{5}\diag(5,4,4,5,9,8,9)$, $||\beta||^2=\frac{7}{5}\thickapprox1.400$
  \newline $[{e_1},{e_2}]=\frac{\sqrt {5}}{5}{e_5},[{e_2},{e_3}]=\frac{\sqrt {10}}{10}{e_6},[{e_3},{e_4}]=\frac{\sqrt {5}}{5}{e_7}$

  \item \textbf{(37B1)}: $\dim(\Der)=20$, $\rank=4$, Betti Numbers $(4, 11, 16, 16, 11, 4, 1)$ ($\cong_{\CC} \ggo_{4.1}$).
  \newline Einstein der. $\phi=\frac{1}{5}\diag(4,5,5,4,9,9,8)$, $||\beta||^2=\frac{7}{5}\thickapprox1.400$
  \newline $[{e_1},{e_2}]=\frac{\sqrt{10}}{10}{e_5},[{e_1},{e_3}]=\frac{\sqrt{10}}{10}{e_6},[{e_1},{e_4}]=\frac{\sqrt{10}}{10}{e_7},[{e_2},{e_4}]=\frac{\sqrt{10}}{10}{e_6},[{e_3},{e_4}]=-\frac{\sqrt{10}}{10}{e_5} $

\item \textbf{(37C)}: $\dim(\Der)=22$, $\rank=3$, Betti Numbers $(4, 11, 17, 17, 11, 4, 1)$ ($\cong_{\RR} \ggo_{3.24}$).
  \newline Einstein der. $\phi=\frac{1}{4}\diag(5,3,4,4,8,7,7)$, $||\beta||^2=\frac{3}{2}\thickapprox1.500$
  \newline $[{e_1},{e_2}]=\frac{\sqrt{2}}{4}{e_5},[{e_2},{e_3}]=\frac{\sqrt{2}}{4}{e_6},[{e_2},{e_4}]=\frac{\sqrt{2}}{4}{e_7},[{e_3},{e_4}]=\frac{\sqrt{2}}{4}{e_5} $

\item \textbf{(37D)}: $\dim(\Der)=19$, $\rank=3$, Betti Numbers $(4, 11, 14, 14, 11, 4, 1)$ ($\cong_{\RR} \ggo_{3.12}$).
  \newline Einstein der. $\phi=\frac{5}{6}\diag(1,1,1,1,2,2,2)$, $||\beta||^2=\frac{4}{3}\thickapprox 1.333$
  \newline $ [{e_1},{e_2}]=\frac{\sqrt{3}}{6}{e_5},[{e_1},{e_3}]=\frac{\sqrt{6}}{6}{e_6},[{e_2},{e_4}]=\frac{\sqrt{6}}{6}{e_7},[{e_3},{e_4}]=\frac{\sqrt{3}}{6}{e_5}$

\item \textbf{(37D1)}: $\dim(\Der)=19$, $\rank=3$, Betti Numbers $(4, 11, 14, 14, 11, 4, 1)$ ($\cong_{\CC} \ggo_{3.12}$).
  \newline Einstein der. $\phi=\frac{5}{6}\diag(1,1,1,1,2,2,2)$, $||\beta||^2=\frac{4}{3}\thickapprox 1.333$
  \newline $[{e_1},{e_2}]=\frac{\sqrt{3}}{6}{e_5},[{e_1},{e_3}]=\frac{\sqrt{3}}{6}{e_6},[{e_1},{e_4}]=\frac{\sqrt{3}}{6}{e_7},[{e_2},{e_3}]=-\frac{\sqrt{3}}{6}{e_7},[{e_2},{e_4}]=\frac{\sqrt{3}}{6}{e_6},[{e_3},{e_4}]=-\frac{\sqrt{3}}{6}{e_5} $


\item \textbf{(357A)}: $\dim(\Der)=18$, $\rank=3$, Betti Numbers $(3, 8, 14, 14, 8, 3, 1)$ ($\cong_{\RR} \ggo_{3.6}$).
  \newline Einstein der. $\phi=\frac{1}{11}\diag(5,7,12,9,17,16,14)$, $||\beta||^2=\frac{13}{11}\thickapprox 1.182$
  \newline $ [{e_1},{e_2}]=\frac{\sqrt {66}}{22}{e_3},[{e_1},{e_3}]=\frac{\sqrt {22}}{11}{e_5},[{e_1},{e_4}]=\frac{\sqrt {22}}{22}{e_7},[{e_2},{e_4}]=\frac{\sqrt {66}}{22}{e_6}$

\item \textbf{(357B)}: $\dim(\Der)=17$, $\rank=3$, Betti Numbers $(3, 7, 11, 11, 7, 3, 1)$ ($\cong_{\RR} \ggo_{3.23}$).
  \newline Einstein der. $\phi=\frac{1}{10}\diag(4,5,9,9,13,14,13)$, $||\beta||^2=\frac{11}{10}\thickapprox1.100$
  \newline $[{e_1},{e_2}]=\frac{\sqrt {15}}{10}{e_3},[{e_1},{e_3}]=\frac{\sqrt {10}}{10}{e_5},[{e_1},{e_4}]=\frac{\sqrt {10}}{10}{e_7},[{e_2},{e_3}]=\frac{\sqrt {15}}{10}{e_6}$

\item \textbf{(357C)}: $\dim(\Der)=16$, $\rank=2$, Betti Numbers $(3, 7, 11, 11, 7, 3, 1)$ ($\cong_{\RR} \ggo_{2.40}$).
  \newline Einstein der. $\phi=\frac{1}{21}\diag(9,10,19,18,28,29,27)$, $||\beta||^2=\frac{23}{21}\thickapprox  1.095$
  \newline $ [{e_1},{e_2}]=\frac{\sqrt{7}}{7}{e_3},[{e_1},{e_3}]=\frac{\sqrt{42}}{21}{e_5},[{e_1},{e_4}]=\frac{\sqrt{42}}{21}{e_7},[{e_2},{e_3}]=\frac{\sqrt{7}}{7}{e_6},
  [{e_2},{e_4}]=\frac{\sqrt {42}}{42}{e_5}$


\item \textbf{(27A)}: $\dim(\Der)=21$, $\rank=4$, Betti Numbers $(5, 10, 16, 16, 10, 5, 1)$ ($\cong_{\RR} \ggo_{4.3}$).
  \newline Einstein der. $\phi=\frac{1}{5}\diag(4,5,5,6,5,9,10)$, $||\beta||^2=\frac{7}{5}\thickapprox 1.400$
  \newline $[{e_1},{e_2}]=\frac{\sqrt{5}}{5}{e_6},[{e_1},{e_4}]=\frac{\sqrt{10}}{10}{e_7},[{e_3},{e_5}]=\frac{\sqrt{5}}{5}{e_7} $

\item \textbf{(27B)}: $\dim(\Der)=19$, $\rank=3$, Betti Numbers $(5, 9, 15, 15, 9, 5, 1)$ ($\cong_{\RR} \ggo_{3.19}$).
  \newline Einstein der. $\phi=\frac{1}{6}\diag(5,6,5,6,6,11,11)$, $||\beta||^2=\frac{4}{3}\thickapprox 1.333$
  \newline $[{e_1},{e_2}]=\frac{\sqrt{3}}{6}{e_6},[{e_1},{e_5}]=\frac{\sqrt{6}}{6}{e_7},[{e_2},{e_3}]=\frac{\sqrt{3}}{6}{e_7},[{e_3},{e_4}]=\frac{\sqrt{6}}{6}{e_6}$


\item \textbf{(257A)}: $\dim(\Der)=19$, $\rank=3$, Betti Numbers $(4, 9, 14, 14, 9, 4, 1)$ ($\cong_{\RR} \ggo_{3.8}$).
  \newline Einstein der. $\phi=\frac{1}{4}\diag(2,3,5,4,4,7,6)$, $||\beta||^2=\frac{5}{4}\thickapprox  1.250 $
  \newline $[{e_1},{e_2}]=\frac{\sqrt{2}}{4}{e_3},[{e_1},{e_3}]=\frac{\sqrt{2}}{4}{e_6},[{e_1},{e_5}]=\frac{\sqrt{2}}{4}{e_7},[{e_2},{e_4}]=\frac{\sqrt{2}}{4}{e_6}$

\item \textbf{(257B)}: $\dim(\Der)=18$, $\rank=3$, Betti Numbers $(4, 8, 13, 13, 8, 4, 1)$ ($\cong_{\RR} \ggo_{3.11}$).
  \newline Einstein der. $\phi=\frac{1}{11}\diag(5,7,12,12,10,17,17)$, $||\beta||^2=\frac{13}{11}\thickapprox  1.182$
  \newline $
  [{e_1},{e_2}]=\frac{\sqrt{66}}{22}{e_3},[{e_1},{e_3}]=\frac{\sqrt{22}}{11}{e_6},[{e_1},{e_4}]=\frac{\sqrt{22}}{22}{e_7},[{e_2},{e_5}]=\frac{\sqrt{66}}{22}{e_7}$

\item \textbf{(257C)}: $\dim(\Der)=18$, $\rank=3$, Betti Numbers $(4, 9, 13, 13, 9, 4, 1)$ ($\cong_{\RR} \ggo_{3.9}$).
  \newline Einstein der. $\phi=\frac{2}{11}\diag(3,3,6,6,5,9,8)$, $||\beta||^2=\frac{13}{11}\thickapprox  1.182$
  \newline $ [{e_1},{e_2}]=\frac{\sqrt{66}}{22}{e_3},[{e_1},{e_3}]=\frac{\sqrt{22}}{11}{e_6},[{e_2},{e_4}]=\frac{\sqrt{22}}{22}{e_6},[{e_2},{e_5}]=\frac{\sqrt{66}}{22}{e_7}$
\item \textbf{(257D)}: $\dim(\Der)=17$, $\rank=2$, Betti Numbers $(4, 8, 12, 12, 8, 4, 1)$ ($\cong_{\RR} \ggo_{2.45}$).
  \newline Einstein der. $\phi=\frac{1}{12}\diag(6,7,13,12,11,19,18)$, $||\beta||^2=\frac{7}{6}\thickapprox 1.167$
  \newline $ [{e_1},{e_2}]=\frac{\sqrt{2}}{4}{e_3},[{e_1},{e_3}]=\frac{\sqrt{6}}{6}{e_6},[{e_1},{e_4}]=\frac{\sqrt{6}}{12}{e_7},[{e_2},{e_4}]=\frac{\sqrt{6}}{12}{e_6},[{e_2},{e_5}]=\frac{\sqrt{2}}{4}{e_7}
$

\item \textbf{(257E)}: $\dim(\Der)=17$, $\rank=3$, Betti Numbers $(4, 8, 11, 11, 8, 4, 1)$ ($\cong_{\RR} \ggo_{3.15}$).
  \newline Einstein der. $\phi=\frac{1}{10}\diag(5,6,11,7,9,16,13)$, $||\beta||^2=\frac{11}{10}\thickapprox 1.100$
  \newline $ [{e_1},{e_2}]=\frac{\sqrt{15}}{10}{e_3},[{e_1},{e_3}]=\frac{\sqrt{15}}{10}{e_6},[{e_2},{e_4}]=\frac{\sqrt{10}}{10}{e_7},[{e_4},{e_5}]=\frac{\sqrt{10}}{10}{e_6}
$

\item \textbf{(257F)}: $\dim(\Der)=18$, $\rank=3$, Betti Numbers $(4, 9, 12, 12, 9, 4, 1)$ ($\cong_{\RR} \ggo_{3.14}$).
  \newline Einstein der. $\phi=\frac{1}{11}\diag(9,5,14,9,10,19,14)$, $||\beta||^2=\frac{13}{11}\thickapprox 1.182$
  \newline $[{e_1},{e_2}]=\frac{\sqrt{22}}{11}{e_3},[{e_2},{e_3}]=\frac{\sqrt{66}}{22}{e_6},[{e_2},{e_4}]=\frac{\sqrt{22}}{22}{e_7},[{e_4},{e_5}]=\frac{\sqrt{66}}{22}{e_6}$

\item \textbf{(257G)}: $\dim(\Der)=16$, $\rank=2$, Betti Numbers $(4, 8, 11, 11, 8, 4, 1)$ ($\cong_{\RR} \ggo_{2.36}$).
  \newline Einstein der. $\phi=\frac{1}{21}\diag(10,13,23,15,18,33,28)$, $||\beta||^2=\frac{23}{21}\thickapprox 1.095$
  \newline $[{e_1},{e_2}]=\frac{\sqrt{7}}{7}{e_3},[{e_1},{e_3}]=\frac{\sqrt{7}}{7}{e_6},[{e_1},{e_5}]=\frac{\sqrt{42}}{42}{e_7},[{e_2},{e_4}]=\frac{\sqrt{42}}{21}{e_7},[{e_4},{e_5}]=\frac{\sqrt{42}}{21}{e_6}
$

\item \textbf{(257H)}: $\dim(\Der)=15$, $\rank=3$, Betti Numbers $(4, 8, 11, 11, 8, 4, 1)$ ($\cong_{\RR} \ggo_{3.7}$).
  \newline Pre-Einstein der. $\phi=\frac{1}{3}\diag(1, 2, 3, 2, 2, 4, 4)$
  \newline It does not admit nilsoliton metrics.

\item \textbf{(257I)!}: $\dim(\Der)=17$, $\rank=2$, Betti Numbers $(4, 8, 11, 11, 8, 4, 1)$ ($\cong_{\RR} \ggo_{2.27}$!).
  \newline Einstein der. $\phi=\frac{1}{13}\diag(6,7,13,13,14,19,20)$, $||\beta||^2=\frac{15}{13}\thickapprox1.154$
  \newline $[{e_1},{e_2}]=\frac{\sqrt{78}}{26}{e_3}+\frac{\sqrt{26}}{26}{e_4},[{e_1},{e_3}]=\frac{\sqrt{26}}{26}{e_6},[{e_1},{e_4}]=\frac{\sqrt{78}}{26}{e_6},[{e_1},{e_5}]=\frac{\sqrt{26}}{26}{e_7},[{e_2},{e_3}]=\frac{\sqrt{26}}{13}{e_7}$

\item \textbf{(257J)}: $\dim(\Der)=16$, $\rank=2$, Betti Numbers $(4, 8, 11, 11, 8, 4, 1)$ ($\cong_{\RR} \ggo_{2.38}$).
  \newline Einstein der. $\phi=\frac{1}{2}\diag(1,1,2,2,2,3,3)$, $||\beta||^2=\frac{8}{7}\thickapprox 1.143$
  \newline
$[{e_1},{e_2}]=\frac{\sqrt{7}}{7}{e_3},[{e_1},{e_3}]=\frac{\sqrt{21}}{14}{e_6},[{e_1},{e_5}]=\frac{\sqrt{14}}{14}{e_7},[{e_2},{e_3}]=\frac{\sqrt{21}}{14}{e_7},[{e_2},{e_4}]=\frac{\sqrt{14}}{14}{e_6}$

\item \textbf{(257J1)!}: $\dim(\Der)=16$, $\rank=2$, Betti Numbers $(4, 8, 11, 11, 8, 4, 1)$ ($\cong_{\CC} \ggo_{2.38}$).
  \newline Einstein der. $\phi=\frac{1}{2}\diag(1,1,2,2,2,3,3)$, $||\beta||^2=\frac{8}{7}\thickapprox 1.143$
  \newline $[{e_1},{e_2}]=\frac{\sqrt{21}}{14}{e_3}+\frac{\sqrt{7}}{14}{e_4},[{e_1},{e_3}]=\frac{\sqrt{7}}{14}{e_6},[{e_1},{e_4}]=\frac{\sqrt{21}}{14}{e_6},[{e_1},{e_5}]=\frac{\sqrt{7}}{14}{e_7},[{e_2},{e_3}]=\frac{\sqrt{7}}{7}{e_7},[{e_2},{e_5}]=\frac{\sqrt{7}}{14}{e_6}
$

\item \textbf{(257K)}: $\dim(\Der)=16$, $\rank=3$, Betti Numbers $(4, 6, 9, 9, 6, 4, 1)$ ($\cong_{\RR} \ggo_{3.13}$).
  \newline Einstein der. $\phi=\frac{1}{22}\diag(8,11,19,15,15,27,30)$, $||\beta||^2=\frac{21}{22}\thickapprox 0.9546$
  \newline $ [{e_1},{e_2}]=\frac{\sqrt{77}}{22}{e_3},[{e_1},{e_3}]=\frac{\sqrt{66}}{22}{e_6},[{e_2},{e_3}]=\frac{\sqrt{33}}{22}{e_7},[{e_4},{e_5}]=\frac{\sqrt{66}}{22}{e_7}$

\item \textbf{(257L)}: $\dim(\Der)=14$, $\rank=2$, Betti Numbers $(4, 6, 9, 9, 6, 4, 1)$ ($\cong_{\RR} \ggo_{2.29}$).
  \newline Pre-Einstein der. $\phi=\frac{1}{41}\diag(15, 22, 37, 30, 29, 52, 59)$,
  \newline It does not admit.


\item \textbf{(247A)}: $\dim(\Der)=19$, $\rank=3$, Betti Numbers $(3, 7, 13, 13, 7, 3, 1)$ ($\cong_{\RR} \ggo_{3.20}$).
  \newline Einstein der. $\phi=\frac{1}{4}\diag(1,4,4,5,5,6,6)$, $||\beta||^2=\frac{5}{4}\thickapprox1.250$
  \newline $[{e_1},{e_2}]=\frac{\sqrt{2}}{4}{e_4},[{e_1},{e_3}]=\frac{\sqrt{2}}{4}{e_5},[{e_1},{e_4}]=\frac{\sqrt{2}}{4}{e_6},[{e_1},{e_5}]=\frac{\sqrt{2}}{4}{e_7}$

\item \textbf{(247B)}: $\dim(\Der)=15$, $\rank=3$, Betti Numbers $(3, 6, 10, 10, 6, 3, 1)$ ($\cong_{\RR} \ggo_{3.21}$).
  \newline Einstein der. $\phi=\frac{1}{22}\diag(6,15,11,21,17,27,28)$, $||\beta||^2=\frac{21}{22}\thickapprox0.9546$
  \newline $ [{e_1},{e_2}]=\frac{\sqrt{66}}{22}{e_4},[{e_1},{e_3}]=\frac{\sqrt{33}}{22}{e_5},[{e_1},{e_4}]=\frac{\sqrt{66}}{22}{e_6},[{e_3},{e_5}]=\frac{\sqrt{77}}{22}{e_7}
$

\item \textbf{(247C)}: $\dim(\Der)=16$, $\rank=2$, Betti Numbers $(3, 7, 11, 11, 7, 3, 1)$ ($\cong_{\RR} \ggo_{2.43}$).
  \newline Einstein der. $\phi=\frac{1}{35}\diag(11,29,20,40,31,51,42)$, $||\beta||^2=\frac{37}{35}\thickapprox1.057$
  \newline $[{e_1},{e_2}]={\frac {2\sqrt {35}}{35}}{e_4},[{e_1},{e_3}]={\frac {2\sqrt {35}}{35}}{e_5},[{e_1},{e_4}]=\frac{\sqrt{14}}{14}{e_6},[{e_1},{e_5}]=\frac{\sqrt{14}}{14}{e_7},[{e_3},{e_5}]={\frac {3\sqrt {70}}{70}}{e_6} $

\item \textbf{(247D)}: $\dim(\Der)=15$, $\rank=3$, Betti Numbers $(3, 6, 10, 10, 6, 3, 1)$ ($\cong_{\RR} \ggo_{3.22}$).
  \newline Einstein der. $\phi=\frac{1}{22}\diag(7,10,12,17,19,24,29)$, $||\beta||^2=\frac{10}{11}\thickapprox 0.9091$
  \newline $ [{e_1},{e_2}]=\frac{\sqrt{55}}{22}{e_4},[{e_1},{e_3}]=\frac{\sqrt{11}}{11}{e_5},[{e_1},{e_4}]=\frac{\sqrt{11}}{11}{e_6},[{e_2},{e_5}]=\frac{\sqrt{55}}{22}{e_7},[{e_3},{e_4}]=\frac{\sqrt{11}}{11}{e_7}$

\item \textbf{(247E)!}: $\dim(\Der)=14$, $\rank=2$, Betti Numbers $(3, 5, 9, 9, 5, 3, 1)$ ($\cong_{\CC} \ggo_{2.12}$).
  \newline Einstein der. $\phi=\frac{1}{10}\diag(3,5,5,8,8,11,13)$, $||\beta||^2=\frac{9}{10}\thickapprox 0.9000$
  \newline $ [{e_1},{e_2}]=\frac{\sqrt{2}}{4}{e_4},[{e_1},{e_3}]=\frac{\sqrt{30}}{20}{e_5},[{e_1},{e_4}]=\frac{\sqrt{10}}{10}{e_6},[{e_2},{e_4}]=\frac{\sqrt{30}}{20}{e_7},[{e_3},{e_5}]=-\frac{\sqrt{2}}{4}{e_7}$

\item \textbf{(247E1)}: $\dim(\Der)=14$, $\rank=2$, Betti Numbers $(3, 5, 9, 9, 5, 3, 1)$ ($\cong_{\RR} \ggo_{2.12}$).
  \newline Einstein der. $\phi=\frac{1}{10}\diag(3,5,5,8,8,11,13)$, $||\beta||^2=\frac{9}{10}\thickapprox0.9000$
  \newline $[{e_1},{e_2}]=\frac{\sqrt{2}}{4}{e_4},[{e_1},{e_3}]=\frac{\sqrt{30}}{20}{e_5},[{e_1},{e_4}]=\frac{\sqrt{10}}{10}{e_6},[{e_2},{e_4}]=\frac{\sqrt{30}}{20}{e_7},[{e_3},{e_5}]=\frac{\sqrt{2}}{4}{e_7}$

\item \textbf{(247F)}: $\dim(\Der)=13$, $\rank=3$, Betti Numbers $(3, 6, 10, 10, 6, 3, 1)$ ($\cong_{\RR} \ggo_{3.4}$).
  \newline Einstein der. $\phi=\frac{1}{14}\diag(6,5,5,11,11,16,16)$, $||\beta||^2=\frac{6}{7}\thickapprox0.8571$
  \newline $[{e_1},{e_2}]=\frac{\sqrt{21}}{14}{e_4},[{e_1},{e_3}]=\frac{\sqrt{21}}{14}{e_5},[{e_2},{e_4}]=\frac{\sqrt{14}}{14}{e_6},[{e_2},{e_5}]=\frac{\sqrt{14}}{14}{e_7},[{e_3},{e_4}]=\frac{\sqrt{14}}{14}{e_7},[{e_3},{e_5}]=\frac{\sqrt{14}}{14}{e_6}$

\item \textbf{(247F1)}: $\dim(\Der)=13$, $\rank=3$, Betti Numbers $(3, 6, 10, 10, 6, 3, 1)$ ($\cong_{\CC} \ggo_{3.4}$).
  \newline Einstein der. $\phi=\frac{1}{14}\diag(6,5,5,11,11,16,16)$, $||\beta||^2=\frac{6}{7}\thickapprox0.8571$
  \newline $[{e_1},{e_2}]=\frac{\sqrt{21}}{14}{e_4},[{e_1},{e_3}]=\frac{\sqrt{21}}{14}{e_5},[{e_2},{e_4}]=\frac{\sqrt{14}}{14}{e_6},[{e_2},{e_5}]=\frac{\sqrt{14}}{14}{e_7},[{e_3},{e_4}]=\frac{\sqrt{14}}{14}{e_7},[{e_3},{e_5}]=-\frac{\sqrt{14}}{14}{e_6} $

\item \textbf{(247G)!}: $\dim(\Der)=12$, $\rank=2$, Betti Numbers $(3, 5, 9, 9, 5, 3, 1)$ ($\cong_{\RR} \ggo_{2.34}$).
  \newline Einstein der. $\phi=\frac{1}{55}\diag(22,20,21,42,43,62,64)$, $||\beta||^2=\frac{47}{55}\thickapprox 0.8546$
  \newline $[{e_1},{e_2}]=-{\frac {\sqrt {330}}{55}}{e_4},[{e_1},{e_3}]=-\frac{\sqrt{10}}{10}{e_5},[{e_1},{e_4}]={\frac {\sqrt {55}}{55}}{e_7},[{e_2},{e_4}]=\frac{\sqrt{66}}{22}{e_6},[{e_3},{e_5}]=-\frac{\sqrt{66}}{22}{e_7}$

\item \textbf{(247H)!}: $\dim(\Der)=11$, $\rank=1$, Betti Numbers $(3, 5, 9, 9, 5, 3, 1)$ ($\cong_{\RR} \ggo_{1.19}$).
  \newline Einstein der. $\phi=\frac{13}{34}\diag(1,1,1,2,2,3,3)$, $||\beta||^2=\frac{29}{34}\thickapprox 0.8529$
  \newline $[{e_1},{e_2}]=-\frac{\sqrt{119}}{34}{e_4},[{e_1},{e_3}]=-\frac{\sqrt{119}}{34}{e_5},[{e_1},{e_4}]=\frac{\sqrt{17}}{34}{e_7},[{e_1},{e_5}]=-\frac{\sqrt{17}}{34}{e_6},[{e_2},{e_4}]={\frac {3\sqrt {17}}{34}}{e_6},[{e_3},{e_5}]=-{\frac {3\sqrt {17}}{34}}{e_7}
$

\item \textbf{(247H1)!}: $\dim(\Der)=11$, $\rank=1$, Betti Numbers $(3, 5, 9, 9, 5, 3, 1)$ ($\cong_{\CC} \ggo_{1.19}$).
  \newline Einstein der. $\phi=\frac{13}{34}\diag(1,1,1,2,2,3,3)$, $||\beta||^2=\frac{29}{34}\thickapprox 0.8529$
  \newline $ [{e_1},{e_2}]=-\frac{\sqrt{119}}{34}{e_4},[{e_1},{e_3}]={\frac {3\sqrt {1309}}{374}}{e_5},[{e_1},{e_4}]={\frac {3\sqrt {187}}{187}}{e_6},[{e_2},{e_3}]={\frac {\sqrt {2618}}{374}}{e_5},[{e_2},{e_4}]=-{\frac {7\sqrt {374}}{748}}{e_6},[{e_2},{e_5}]=-{\frac {\sqrt {374}}{68}}{e_7},[{e_3},{e_4}]={\frac {3\sqrt {34}}{68}}{e_7},[{e_3},{e_5}]=-{\frac {3\sqrt {34}}{68}}{e_6}$

\item \textbf{(247I)}: $\dim(\Der)=14$, $\rank=3$, Betti Numbers $(3, 6, 10, 10, 6, 3, 1)$ ($\cong_{\RR} \ggo_{3.5}$).
  \newline Einstein der. $\phi=\frac{1}{22}\diag(10,11,7,21,17,28,24)$, $||\beta||^2=\frac{10}{11}\thickapprox0.9091$
  \newline $[{e_1},{e_2}]=\frac{\sqrt{55}}{22}{e_4},[{e_1},{e_3}]=\frac{\sqrt{55}}{22}{e_5},[{e_2},{e_5}]=\frac{\sqrt{11}}{11}{e_6},[{e_3},{e_4}]=\frac{\sqrt{11}}{11}{e_6},[{e_3},{e_5}]=\frac{\sqrt{11}}{11}{e_7}$

\item \textbf{(247J)!}: $\dim(\Der)=13$, $\rank=2$, Betti Numbers $(3, 6, 10, 10, 6, 3, 1)$ ($\cong_{\RR} \ggo_{2.26}$!).
  \newline Einstein der. $\phi=\frac{1}{19}\diag(7,10,7,17,14,21,24)$, $||\beta||^2=\frac{17}{19}\thickapprox 0.8947$
  \newline $[{e_1},{e_2}]=\frac{\sqrt{114}}{38}{e_4},[{e_1},{e_3}]=\frac{\sqrt{38}}{19}{e_5},[{e_1},{e_5}]=\frac{\sqrt{114}}{38}{e_6},[{e_2},{e_3}]=\frac{\sqrt{38}}{38}{e_4},[{e_2},{e_5}]=\frac{\sqrt{114}}{38}{e_7},[{e_3},{e_4}]=\frac{\sqrt{38}}{19}{e_7},[{e_3},{e_5}]=\frac{\sqrt{38}}{38}{e_6}
$

\item \textbf{(247K)}: $\dim(\Der)=12$, $\rank=2$, Betti Numbers $(3, 5, 9, 9, 5, 3, 1)$ ($\cong_{\RR} \ggo_{2.35}$).
  \newline Einstein der. $\phi=\frac{1}{15}\diag(5,7,6,12,11,17,18)$, $||\beta||^2=\frac{13}{15}\thickapprox0.8667$
  \newline $[{e_1},{e_2}]={\frac {\sqrt{330}}{55}}{e_4},[{e_1},{e_3}]=\frac{\sqrt{11}}{11}{e_5},[{e_1},{e_4}]=\frac{\sqrt{15}}{15}{e_6},[{e_2},{e_5}]=\frac{\sqrt{11}}{11}{e_7},[{e_3},{e_4}]={\frac {\sqrt{330}}{66}}{e_7},[{e_3},{e_5}]=\frac{\sqrt{15}}{15}{e_6}$

\item \textbf{(247L)}: $\dim(\Der)=17$, $\rank=2$, Betti Numbers $(3, 7, 13, 13, 7, 3, 1)$ ($\cong_{\RR} \ggo_{2.39}$).
  \newline Einstein der. $\phi=\frac{1}{14}\diag(5,11,10,16,15,21,20)$, $||\beta||^2=\frac{8}{7}\thickapprox  1.143$
  \newline $[{e_1},{e_2}]=\frac{\sqrt{14}}{14}{e_4},[{e_1},{e_3}]=\frac{\sqrt{21}}{14}{e_5},[{e_1},{e_4}]=\frac{\sqrt{14}}{14}{e_6},[{e_1},{e_5}]=\frac{\sqrt{7}}{7}{e_7},[{e_2},{e_3}]=\frac{\sqrt{21}}{14}{e_6}$

\item \textbf{(247M)}: $\dim(\Der)=14$, $\rank=2$, Betti Numbers $(3, 6, 10, 10, 6, 3, 1)$ ($\cong_{\RR} \ggo_{2.42}$).
  \newline Pre-Einstein der. $\phi=\frac{1}{41}\diag(11, 30, 22, 41, 33, 52, 55)$.
  \newline It does not admit.

\item \textbf{(247N)}: $\dim(\Der)=16$, $\rank=2$, Betti Numbers $(3, 7, 11, 11, 7, 3, 1)$ ($\cong_{\RR} \ggo_{2.44}$).
  \newline Einstein der. $\phi=\frac{1}{35}\diag(15,19,23,34,38,53,42)$, $||\beta||^2=\frac{37}{35}\thickapprox1.057$
  \newline $[{e_1},{e_2}]=\frac{\sqrt{14}}{14}{e_4},[{e_1},{e_3}]={\frac {3\sqrt{70}}{70}}{e_5},[{e_1},{e_5}]={\frac {2\sqrt{35}}{35}}{e_6},[{e_2},{e_3}]=\frac{\sqrt{14}}{14}{e_7},[{e_2},{e_4}]={\frac {2\sqrt{35}}{35}}{e_6}$

\item \textbf{(247O)}: $\dim(\Der)=15$, $\rank=1$, Betti Numbers $(3, 7, 11, 11, 7, 3, 1)$ ($\cong_{\RR} \ggo_{1.7}$).
  \newline Einstein der. $\phi=\frac{5}{28}\diag(2,4,3,6,5,8,7)$, $||\beta||^2=\frac{29}{28}\thickapprox 1.036$
  \newline $ [{e_1},{e_2}]=\frac{\sqrt{21}}{14}{e_4},[{e_1},{e_3}]=\frac{\sqrt{70}}{28}{e_5},[{e_1},{e_4}]=\frac{\sqrt{70}}{28}{e_6},[{e_1},{e_5}]=\frac{\sqrt{42}}{28}{e_7},[{e_2},{e_3}]=\frac{\sqrt{42}}{28}{e_7},[{e_3},{e_5}]=\frac{\sqrt{21}}{14}{e_6}$

\item \textbf{(247P)}: $\dim(\Der)=15$, $\rank=3$, Betti Numbers $(3, 7, 11, 11, 7, 3, 1)$ ($\cong_{\RR} \ggo_{3.1(i_{0})}$).
  \newline Pre-Einstein der. $\phi=\frac{1}{2}\diag(1, 1, 1, 2, 2, 2, 3)$,
  \newline It does not admit.

\item \textbf{(247P1)}: $\dim(\Der)=15$, $\rank=3$, Betti Numbers $(3, 7, 11, 11, 7, 3, 1)$ ($\cong_{\CC} \ggo_{3.1(i_{0})}$).
  \newline Pre-Einstein der. $\phi=\frac{1}{2}\diag(1, 1, 1, 2, 2, 2, 3)$,
  \newline It does not admit.

\item \textbf{(247Q)}: $\dim(\Der)=14$, $\rank=2$, Betti Numbers $(3, 6, 10, 10, 6, 3, 1)$ ($\cong_{\RR} \ggo_{2.1(v)}$).
  \newline Pre-Einstein der. $\phi=\frac{2}{19}\diag(3, 5, 6, 8, 9, 11, 14)$,
  \newline It does not admit.

\item \textbf{(247R)!}: $\dim(\Der)=13$, $\rank=1$, Betti Numbers $(3, 5, 9, 9, 5, 3, 1)$ ($\cong_{\CC} \ggo_{1.3(iv)}$).
  \newline Pre-Einstein der. $\phi=\frac{5}{17}\diag(1, 2, 2, 3, 3, 4, 5)$,
  \newline It does not admit.

\item \textbf{(247R1)}: $\dim(\Der)=13$, $\rank=1$, Betti Numbers $(3, 5, 9, 9, 5, 3, 1)$ ($\cong_{\RR} \ggo_{1.3(iv)}$).
  \newline Pre-Einstein der. $\phi=\frac{5}{17}\diag(1, 2, 2, 3, 3, 4, 5)$,
  \newline It does not admit


\item \textbf{(2457A)}: $\dim(\Der)=17$, $\rank=3$, Betti Numbers $(3, 7, 10, 10, 7, 3, 1)$ ($\cong_{\RR} \ggo_{3.2}$).
  \newline Einstein der. $\phi=\frac{2}{11}\diag(1,5,6,7,6,8,7)$, $||\beta||^2=\frac{13}{11}\thickapprox  1.182$
  \newline $[{e_1},{e_2}]=\frac{\sqrt{66}}{22}{e_3},[{e_1},{e_3}]=\frac{\sqrt{22}}{11}{e_4},[{e_1},{e_4}]=\frac{\sqrt{66}}{22}{e_6},[{e_1},{e_5}]=\frac{\sqrt{22}}{22}{e_7}$

\item \textbf{(2457B)}: $\dim(\Der)=15$, $\rank=3$, Betti Numbers $(3, 7, 9, 9, 7, 3, 1)$ ($\cong_{\RR} \ggo_{3.3}$).
  \newline Einstein der. $\phi=\frac{1}{22}\diag(5,12,17,22,15,27,27)$, $||\beta||^2=\frac{21}{22}\thickapprox0.9546$
  \newline $[{e_1},{e_2}]=\frac{\sqrt{33}}{22}{e_3},[{e_1},{e_3}]=\frac{\sqrt{77}}{22}{e_4},[{e_1},{e_4}]=\frac{\sqrt{66}}{22}{e_7},[{e_2},{e_5}]=\frac{\sqrt{66}}{22}{e_6}$

\item \textbf{(2457C)}: $\dim(\Der)=19$, $\rank=2$, Betti Numbers $(3, 7, 10, 10, 7, 3, 1)$ ($\cong_{\RR} \ggo_{2.21}$).
  \newline Einstein der. $\phi=\frac{1}{29}\diag(8,19,27,35,24,43,32)$, $||\beta||^2=\frac{31}{29}\thickapprox 1.069$
  \newline $ [{e_1},{e_2}]=\frac{\sqrt{87}}{29}{e_3},[{e_1},{e_3}]=\frac{\sqrt{145}}{29}{e_4},[{e_1},{e_4}]=\frac{\sqrt{87}}{29}{e_6},[{e_1},{e_5}]=\frac{\sqrt{58}}{58}{e_7},[{e_2},{e_5}]=\frac{\sqrt{87}}{29}{e_6}
$

\item \textbf{(2457D)!}: $\dim(\Der)=15$, $\rank=1$, Betti Numbers $(3, 7, 10, 10, 7, 3, 1)$ ($\cong_{\RR} \ggo_{1.16}$!).
  \newline Einstein der. $\phi=\frac{11}{38}\diag(1,2,3,3,4,4,5)$, $||\beta||^2=\frac{20}{19}\thickapprox  1.053$
  \newline $[{e_1},{e_2}]=\frac{\sqrt{114}}{38}{e_3}+\frac{\sqrt{57}}{38}{e_4},[{e_1},{e_3}]=\frac{\sqrt{190}}{38}{e_6},[{e_1},{e_4}]=\frac{\sqrt{19}}{19}{e_5},[{e_1},{e_6}]=\frac{\sqrt{114}}{38}{e_7},[{e_2},{e_3}]=\frac{\sqrt{57}}{38}{e_7},[{e_2},{e_4}]=\frac{\sqrt{114}}{38}{e_7}
$

\item \textbf{(2457E)!}: $\dim(\Der)=14$, $\rank=2$, Betti Numbers $(3, 6, 8, 8, 6, 3, 1)$ ($\cong_{\RR} \ggo_{2.11}$!).
  \newline Einstein der. $\phi=\frac{1}{38}\diag(9,19,28,28,37,47,46)$, $||\beta||^2={\frac {18}{19}}\thickapprox 0.9474$
  \newline $[{e_1},{e_2}]={\frac {\sqrt {1254}}{209}}{e_3}+{\frac {\sqrt {8778}}{418}}{e_4},[{e_1},{e_3}]={\frac {2\sqrt {1463}}{209}}{e_5},[{e_1},{e_4}]={\frac {3\sqrt {209}}{418}}{e_5},[{e_1},{e_5}]=\frac{\sqrt{190}}{38}{e_7},[{e_2},{e_4}]=\frac{\sqrt{209}}{38}{e_6}$

\item \textbf{(2457F)}: $\dim(\Der)=16$, $\rank=2$, Betti Numbers $(3, 7, 10, 10, 7, 3, 1)$ ($\cong_{\RR} \ggo_{2.20}$).
  \newline Einstein der. $\phi=\frac{2}{35}\diag(5,10,15,20,16,25,21)$, $||\beta||^2=\frac{37}{35}\thickapprox1.057$
  \newline $[{e_1},{e_2}]={\frac {3\sqrt {70}}{70}}{e_3},[{e_1},{e_3}]={\frac {2\sqrt {35}}{35}}{e_4},[{e_1},{e_4}]={\frac {\sqrt {14}}{14}}{e_6},[{e_1},{e_5}]={\frac {\sqrt {14}}{14}}{e_7},[{e_2},{e_3}]={\frac {2\sqrt {35}}{35}}{e_6}$

\item \textbf{(2457G)}: $\dim(\Der)=15$, $\rank=2$, Betti Numbers $(3, 6, 9, 9, 6, 3, 1)$ ($\cong_{\RR} \ggo_{2.19}$).
  \newline Pre-Einstein der. $\phi=\frac{1}{4}\diag(1, 2, 3, 4, 4, 5, 5)$,
  \newline It does not admit.

\item \textbf{(2457H)}: $\dim(\Der)=15$, $\rank=2$, Betti Numbers $(3, 6, 10, 10, 6, 3, 1)$ ($\cong_{\RR} \ggo_{2.18}$).
  \newline Einstein der. $\phi=\frac{1}{70}\diag(20,31,51,71,60,82,91)$, $||\beta||^2=\frac{34}{35}\thickapprox0.9714$
  \newline $[{e_1},{e_2}]={\frac {\sqrt {21}}{14}}{e_3},[{e_1},{e_3}]={\frac {3\sqrt {70}}{70}}{e_4},[{e_1},{e_4}]={\frac {\sqrt {21}}{14}}{e_7},[{e_2},{e_3}]={\frac {\sqrt {10}}{10}}{e_6},[{e_2},{e_5}]={\frac {\sqrt {70}}{35}}{e_7} $

\item \textbf{(2457I)}: $\dim(\Der)=14$, $\rank=2$, Betti Numbers $(3, 7, 9, 9, 7, 3, 1)$ ($\cong_{\RR} \ggo_{2.22}$).
  \newline Einstein der. $\phi=\frac{1}{20}\diag(5,10,15,20,14,25,24)$, $||\beta||^2=\frac{19}{20}\thickapprox0.9500$
  \newline $ [{e_1},{e_2}]=\frac{\sqrt{30}}{20}{e_3},[{e_1},{e_3}]=\frac{\sqrt{15}}{10}{e_4},[{e_1},{e_4}]=\frac{\sqrt{2}}{4}{e_6},[{e_2},{e_3}]=\frac{\sqrt{10}}{20}{e_6},[{e_2},{e_5}]=\frac{\sqrt{2}}{4}{e_7}$

\item \textbf{(2457J)!}: $\dim(\Der)=13$, $\rank=1$, Betti Numbers $(3, 6, 8, 8, 6, 3, 1)$ ($\cong_{\RR} \ggo_{1.18}$!).
  \newline Einstein der. $\phi=\frac{23}{94}\diag(1,2,3,3,4,5,5)$, $||\beta||^2=\frac{89}{94}\thickapprox 0.9468$
  \newline $[{e_1},{e_2}]={\frac {3\sqrt {141}}{188}}{e_3}+{\frac {\sqrt {1551}}{188}}{e_4},[{e_1},{e_3}]={\frac {3\sqrt {517}}{188}}{e_5},[{e_1},{e_4}]={\frac {3\sqrt {47}}{188}}{e_5},[{e_1},{e_5}]=\frac{\sqrt {282}}{47}{e_7},[{e_2},{e_3}]={\frac {\sqrt {94}}{94}}{e_7},[{e_2},{e_4}]={\frac {\sqrt {1222}}{94}}{e_6}$

\item \textbf{(2457K)}: $\dim(\Der)=14$, $\rank=1$, Betti Numbers $(3, 6, 9, 9, 6, 3, 1)$ ($\cong_{\RR} \ggo_{1.9}$).
  \newline Pre-Einstein der. $\phi=\frac{10}{67}\diag(2, 3, 5, 7, 6, 8, 9)$,
  \newline It does not admit.

\item \textbf{(2457L)}: $\dim(\Der)=12$, $\rank=2$, Betti Numbers $(2, 5, 8, 8, 5, 2, 1)$ ($\cong_{\RR} \ggo_{2.9}$).
  \newline Einstein der. $\phi=\frac{9}{34}\diag(1,1,2,3,3,4,4)$, $||\beta||^2=\frac{14}{17}\thickapprox0.8235$
  \newline $[{e_1},{e_2}]=\frac{\sqrt{17}}{17}{e_3},[{e_1},{e_3}]=\frac{\sqrt{119}}{34}{e_4},[{e_1},{e_4}]=-\frac{\sqrt{17}}{17}{e_6},[{e_1},{e_5}]=-\frac{\sqrt{17}}{17}{e_7},[{e_2},{e_3}]=-\frac{\sqrt{119}}{34}{e_5},[{e_2},{e_4}]=\frac{\sqrt{17}}{17}{e_7},[{e_2},{e_5}]=\frac{\sqrt{17}}{17}{e_6}
 $

\item \textbf{(2457L1)}: $\dim(\Der)=12$, $\rank=2$, Betti Numbers $(2, 5, 8, 8, 5, 2, 1)$ ($\cong_{\CC} \ggo_{2.9}$).
  \newline Einstein der. $\phi=\frac{9}{34}\diag(1,1,2,3,3,4,4)$, $||\beta||^2=\frac{14}{17}\thickapprox0.8235$
  \newline $[{e_1},{e_2}]=\frac{\sqrt{17}}{17}{e_3},[{e_1},{e_3}]=\frac{\sqrt{119}}{34}{e_4},[{e_1},{e_4}]=-\frac{\sqrt{17}}{17}{e_6},[{e_1},{e_5}]=-\frac{\sqrt{17}}{17}{e_7},[{e_2},{e_3}]=-\frac{\sqrt{119}}{34}{e_5},[{e_2},{e_4}]=\frac{\sqrt{17}}{17}{e_7},[{e_2},{e_5}]=-\frac{\sqrt{17}}{17}{e_6}$

\item \textbf{(2457M)}: $\dim(\Der)=13$, $\rank=2$, Betti Numbers $(2, 5, 9, 9, 5, 2, 1)$ ($\cong_{\RR} \ggo_{2.8}$).
  \newline Einstein der. $\phi=\frac{1}{14}\diag(3,5,8,11,13,16,14)$, $||\beta||^2=\frac{6}{7}\thickapprox 0.8571$
  \newline $[{e_1},{e_2}]=\frac{\sqrt{14}}{14}{e_3},[{e_1},{e_3}]=\frac{\sqrt{21}}{14}{e_4},[{e_1},{e_4}]=\frac{\sqrt{14}}{14}{e_7},[{e_1},{e_5}]=\frac{\sqrt{14}}{14}{e_6},[{e_2},{e_3}]=\frac{\sqrt{21}}{14}{e_5},[{e_2},{e_4}]=\frac{\sqrt{14}}{14}{e_6}$


\item \textbf{(2357A)!}: $\dim(\Der)=13$, $\rank=2$, Betti Numbers $(3, 6, 8, 8, 6, 3, 1)$ ($\cong_{\RR} \ggo_{2.24}$!).
  \newline Einstein der. $\phi=\frac{1}{19}\diag(5,9,10,14,19,19,24)$, $||\beta||^2=\frac{17}{19}\thickapprox  0.8947$
  \newline $[{e_1},{e_2}]=-\frac{\sqrt{38}}{19}{e_4},[{e_1},{e_4}]=\frac{\sqrt{114}}{38}{e_5}-\frac{\sqrt{38}}{38}{e_6},[{e_1},{e_5}]=\frac{\sqrt{38}}{19}{e_7},[{e_2},{e_3}]=\frac{\sqrt{114}}{38}{e_5}+\frac{\sqrt{38}}{38}{e_6},[{e_3},{e_4}]=\frac{\sqrt{114}}{38}{e_7}$

\item \textbf{(2357B)}: $\dim(\Der)=14$, $\rank=2$, Betti Numbers $(3, 6, 8, 8, 6, 3, 1)$ ($\cong_{\RR} \ggo_{2.1(i_\lambda)}$).
  \newline Pre-Einstein der. $\phi=\frac{2}{19}\diag(3, 5, 6, 8, 11, 9, 14)$,
  \newline It does not admit.

\item \textbf{(2357C)}: $\dim(\Der)=13$, $\rank=2$, Betti Numbers $(3, 6, 7, 7, 6, 3, 1)$ ($\cong_{\RR} \ggo_{2.17}$).
  \newline Einstein der. $\phi=\frac{1}{14}\diag(4,5,8,9,13,14,17)$, $||\beta||^2=\frac{6}{7}\thickapprox0.8571$
  \newline $[{e_1},{e_2}]=\frac{\sqrt{21}}{14}{e_4},[{e_1},{e_4}]=\frac{\sqrt{14}}{14}{e_5},[{e_1},{e_5}]=\frac{\sqrt{21}}{14}{e_7},[{e_2},{e_3}]=\frac{\sqrt{14}}{14}{e_5},[{e_2},{e_4}]=\frac{\sqrt{14}}{14}{e_6},[{e_3},{e_4}]=-\frac{\sqrt{14}}{14}{e_7}$

\item \textbf{(2357D)}: $\dim(\Der)=12$, $\rank=1$, Betti Numbers $(3, 6, 7, 7, 6, 3, 1)$ ($\cong_{\RR} \ggo_{1.2(iii)}$).
  \newline Pre-Einstein der. $\phi=\frac{4}{11}\diag(1, 1, 2, 2, 3, 3, 4)$,
  \newline It does not admit.

\item \textbf{(2357D1)}: $\dim(\Der)=12$, $\rank=1$, Betti Numbers $(3, 6, 7, 7, 6, 3, 1)$ ($\cong_{\CC} \ggo_{1.2(iii)}$).
  \newline Pre-Einstein der. $\phi=\frac{4}{11}\diag(1, 1, 2, 2, 3, 3, 4)$,
  \newline It does not admit.


\item \textbf{(23457A)}: $\dim(\Der)=13$, $\rank=2$, Betti Numbers $(2, 4, 7, 7, 4, 2, 1)$ ($\cong_{\RR} \ggo_{2.7}$).
  \newline Einstein der. $\phi=\frac{1}{20}\diag(3,10,13,16,19,22,23)$, $||\beta||^2=\frac{9}{10}\thickapprox0.9000$
  \newline $[{e_1},{e_2}]=\frac{\sqrt{30}}{20}{e_3},[{e_1},{e_3}]=\frac{\sqrt{30}}{20}{e_4},[{e_1},{e_4}]=\frac{\sqrt{2}}{4}{e_5},[{e_1},{e_5}]=\frac{\sqrt{10}}{10}{e_6},[{e_2},{e_3}]=\frac{\sqrt{2}}{4}{e_7}$

\item \textbf{(23457B)}: $\dim(\Der)=12$, $\rank=2$, Betti Numbers $(2, 3, 4, 4, 3, 2, 1)$ ($\cong_{\RR} \ggo_{2.6}$).
  \newline Einstein der. $\phi=\frac{1}{70}\diag(10,23,33,43,53,76,56)$, $||\beta||^2=\frac{26}{35}\thickapprox0.7429$
  \newline $[{e_1},{e_2}]={\frac {\sqrt {455}}{70}}{e_3},[{e_1},{e_3}]={\frac {2\sqrt {35}}{35}}{e_4},[{e_1},{e_4}]={\frac {\sqrt {455}}{70}}{e_5},[{e_2},{e_3}]=\frac{\sqrt{35}}{35}{e_7},[{e_2},{e_5}]=\frac{\sqrt{105}}{35}{e_6},[{e_3},{e_4}]=-\frac{\sqrt{105}}{35}{e_6}$

\item \textbf{(23457C)}: $\dim(\Der)=12$, $\rank=2$, Betti Numbers $(2, 3, 4, 4, 3, 2, 1)$ ($\cong_{\RR} \ggo_{2.4}$).
  \newline Einstein der. $\phi=\frac{1}{10}\diag(1,4,5,6,7,8,11)$, $||\beta||^2=\frac{26}{35}\thickapprox0.7429$
  \newline $[{e_1},{e_2}]=\frac{\sqrt{105}}{35}{e_3},[{e_1},{e_3}]={\frac {2\sqrt {35}}{35}}{e_4},[{e_1},{e_4}]={\frac {\sqrt {455}}{70}}{e_5},[{e_1},{e_5}]=\frac{\sqrt{35}}{35}{e_6},[{e_2},{e_5}]=\frac{\sqrt{105}}{35}{e_7},[{e_3},{e_4}]=-{\frac {\sqrt {455}}{70}}{e_7}$

\item \textbf{(23457D)}: $\dim(\Der)=11$, $\rank=1$, Betti Numbers $(2, 3, 4, 4, 3, 2, 1)$ ($\cong_{\RR} \ggo_{1.5}$).
  \newline Einstein der. $\phi=\frac{5}{42}\diag(1,3,4,5,6,7,9)$, $||\beta||^2=\frac{31}{42}\thickapprox 0.7381$
  \newline $[{e_1},{e_2}]={\frac {\sqrt {290}}{58}}{e_3},[{e_1},{e_3}]=\frac{\sqrt{21}}{14}{e_4},[{e_1},{e_4}]={\frac {5\sqrt {609}}{406}}{e_5},[{e_1},{e_5}]=\frac{\sqrt{42}}{42}{e_6},[{e_2},{e_3}]=\frac{\sqrt{42}}{42}{e_6},[{e_2},{e_5}]={\frac {\sqrt {609}}{87}}{e_7},[{e_3},{e_4}]=-{\frac {\sqrt {290}}{58}}{e_7}$

\item \textbf{(23457E)!}: $\dim(\Der)=12$, $\rank=1$, Betti Numbers $(2, 4, 7, 7, 4, 2, 1)$ ($\cong_{\RR} \ggo_{1.13}$!).
  \newline Einstein der. $\phi=\frac{13}{68}\diag(1,2,3,4,5,5,6)$, $||\beta||^2=\frac{29}{34}\thickapprox 0.8529$
  \newline $[{e_1},{e_2}]={\frac {3\sqrt {34}}{68}}{e_3},[{e_1},{e_3}]=\frac{\sqrt {119}}{34}{e_4},[{e_1},{e_4}]={\frac {3\sqrt {187}}{187}}{e_5}-{\frac {7\sqrt {374}}{748}}{e_6},[{e_1},{e_5}]={\frac {\sqrt {374}}{68}}{e_7},[{e_2},{e_3}]={\frac {3\sqrt {1309}}{374}}{e_5}+{\frac {\sqrt {2618}}{374}}{e_6},[{e_2},{e_4}]={\frac {3\sqrt {34}}{68}}{e_7}$

\item \textbf{(23457F)!}: $\dim(\Der)=11$, $\rank=1$, Betti Numbers $(2, 3, 4, 4, 3, 2, 1)$ ($\cong_{\RR} \ggo_{1.14}$!).
  \newline Einstein der. $\phi=\frac{9}{58}\diag(1,2,3,4,5,5,7)$, $||\beta||^2={\frac {43}{58}}\thickapprox 0.7414$
  \newline $[{e_1},{e_2}]=-{\frac {\sqrt {519593}}{2378}}{e_3},[{e_1},{e_3}]={\frac {\sqrt {377}}{58}}{e_4},[{e_1},{e_4}]={\frac {3\sqrt {605845438}}{252068}}{e_5}+{\frac {\sqrt {126034}}{6148}}{e_6},[{e_2},{e_3}]={\frac {\sqrt {777722}}{6148}}{e_5}-{\frac {3\sqrt {58406}}{6148}}{e_6},[{e_2},{e_5}]={\frac {\sqrt {126034}}{1189}}{e_7},[{e_3},{e_4}]={\frac {3\sqrt {13079}}{1189}}{e_7}$

\item \textbf{(23457G)}: $\dim(\Der)=10$, $\rank=1$, Betti Numbers $(2, 3, 4, 4, 3, 2, 1)$ ($\cong_{\RR} \ggo_{1.1(iii)}$).
  \newline Einstein der. $\phi=\frac{1}{7}\diag(1,2,3,4,5,6,7)$, $||\beta||^2=\frac{5}{7}\thickapprox0.7143$
  \newline $[{e_1},{e_2}]=\frac{\sqrt{14}}{14}{e_3},[{e_1},{e_3}]=\frac{\sqrt{42}}{21}{e_4},[{e_1},{e_4}]=\frac{\sqrt{14}}{14}{e_5},[{e_1},{e_5}]=\frac{\sqrt{84}}{42}{e_6},[{e_2},{e_3}]=\frac{\sqrt{84}}{42}{e_5},
  [{e_2},{e_4}]=\frac{\sqrt {42}}{42}{e_6},[{e_2},{e_5}]=\frac{\sqrt{14}}{14}{e_7},[{e_3},{e_4}]=-\frac{\sqrt{14}}{14}{e_7}$


\item \textbf{(17)}: $\dim(\Der)=28$, $\rank=4$, Betti Numbers $(6, 14, 14, 14, 14, 6, 1)$ ($\cong_{\RR} \ggo_{4.4}$).
  \newline Einstein der. $\phi=\frac{4}{3}\diag(1,1,1,1,1,1,2)$, $||\beta||^2=\frac{5}{3}\thickapprox 1.667$
  \newline $ [{e_1},{e_2}]=\frac{\sqrt{6}}{6}{e_7},[{e_3},{e_4}]=\frac{\sqrt{6}}{6}{e_7},[{e_5},{e_6}]=\frac{\sqrt{6}}{6}{e_7}$


\item \textbf{(157)}: $\dim(\Der)=19$, $\rank=3$, Betti Numbers $(5, 10, 11, 11, 10, 5, 1)$ ($\cong_{\RR} \ggo_{3.18}$).
  \newline Einstein der. $\phi=\frac{2}{11}\diag(3,4,7,6,5,5,10)$, $||\beta||^2=\frac{13}{11}\thickapprox1.182$
  \newline $[{e_1},{e_2}]=\frac{\sqrt{22}}{11}{e_3},[{e_1},{e_3}]=\frac{\sqrt{66}}{22}{e_7},[{e_2},{e_4}]=\frac{\sqrt{22}}{22}{e_7},[{e_5},{e_6}]=\frac{\sqrt{66}}{22}{e_7}$
\end{enumerate}

\begin{example}\label{ex1}
Consider the one-parameter family $(1357QRS1)[t\in\RR]$ given by $(\RR^7,\mu_{t})$ with
$$
\mu_{t}:=\left\{\begin{array}{l}[e_1, e_2] = e_3, [e_1, e_3] = e_5, [e_1, e_4] = e_6, [e_1, e_5] = e_7, \\{[e_2, e_3]} = -e_6, [e_2, e_4] = e_5,[e_2, e_6] = te_7, [e_3, e_4] = (1 - t)e_7\end{array}\right. .
$$
The basis $\{e_1,\ldots,e_7\}$ is a nice basis to $(1357QRS1)[t\in\RR]$ for any $t$. If $t\neq0,1$, the Gram matrix $U$ is given by
$${\tiny \left[
\begin{array}{cccccccc} 3&0&1&1&0&1&1&-1\\
\noalign{\medskip}0&3&1&0&1&1&0&1\\
\noalign{\medskip}1&1&3&1&1&1&-1&1\\
\noalign{\medskip}1&0&1&3&0&-1&1&1\\
\noalign{\medskip}0&1&1&0&3&1&0&1\\
\noalign{\medskip}1&1&1&-1&1&3&1&1\\
\noalign{\medskip}1&0&-1&1&0&1&3&1\\
\noalign{\medskip}-1&1&1&1&1&1&1&3
\end{array} \right]}.
$$
The general solution to the problem $Ux=[1]_{8}$ is given by
$$
x=\frac{1}{11}(1+11\,a_{{1}}, 2, -1+11\,a_{{2}} , 5-11\,a_{{1}}-11\,a_{{2}} , 2 , 4-11\,a_{{1}}-11\,a_{{2}} , 11\,a_{{2}} , 11\,a_{{1}})^{T}.
$$
By taking $a_{1}$ and $a_{2}$ such that $0<11a_{1}<3$ and $1<11a_{2}<4-11a_{1}$, we get a solution with
positive coordinates. Hence, $(1357QRS1)[t\in\RR]$ with $t\neq0,1$ admits a nilsoliton metric (\cite[Theorem 3]{NIKOLAYEVSKY2}).
As we have mentioned above, to give a nilsoliton metric for each $t$ in a simultaneous manner, it is in general incredibly difficult.
But, if we fix the value of $t$, say $t=-1$, it is easy to give the nilsoliton metric for $(1357QRS1)[t=-1]$. To find such metric,
we solve the problem
\begin{eqnarray*}
\mm(g\cdot\mu_{-1}) &=& \mcc(\mathfrak{R}(\mu_{-1}))  \\
                          &=& \frac{1}{13}\diag(-7,-7,-3,-3,1,1,5)
\end{eqnarray*}
 for $g \in \GrpA$ (Theorem \ref{niceTh}). Let $g=\diag(1,1,{\frac {2\sqrt {39}}{39}},\frac{2\sqrt {13}}{13},{\frac {2\sqrt {3}}{39}},{\frac {2\sqrt {3}}{39}},{\frac {\sqrt {130}}{507}})$, the change of basis given by $g$ defines $\widetilde{\mu}_{-1}:=g\cdot \mu_{-1}=g\mu_{-1}(g^{-1}\cdot,g^{-1}\cdot)$
$$
\widetilde{\mu}_{-1}:=\left\{\begin{array}{l} [{e_1},{e_2}]={\frac {2\sqrt {39}}{39}}{e_3},[{e_1},{e_3}]=\frac{\sqrt{13}}{13}{e_5},[{e_1},{e_4}]=\frac{\sqrt{39}}{39}{e_6},[{e_1},{e_5}]={\frac {\sqrt {390}}{78}}{e_7},\\{[{e_2},{e_3}]}=-\frac{\sqrt{13}}{13}{e_6},[{e_2},{e_4}]=\frac{\sqrt{39}}{39}{e_5},[{e_2},{e_6}]=-{\frac {\sqrt {390}}{78}}{e_7},[{e_3},{e_4}]={\frac {\sqrt {390}}{78}}{e_7}
\end{array}\right. .
$$
Since,
\begin{eqnarray}
\nonumber\label{eq1} \mm(\widetilde{\mu}_{-1}) &=& \frac{1}{13}\diag(-7,-7,-3,-3,1,1,5) \\
                          &=& -\frac{11}{13}\Id + \underbrace{\frac{4}{13}\diag(1,1,2,2,3,3,4)}_{\mbox{Derivation}}
\end{eqnarray}
the canonical inner product of $\RR^7$ defines a nilsoliton metric on $(\RR^7,\widetilde{\mu}_{-1})$. Note that $(1357QRS1)[t=-1]$ becomes $(1357R)$ over $\CC$, therefore $(1357R)$ also admits nilsoliton metrics (Theorem \ref{complex1}).

We study $(1357QRS1)[t=1]$ ($\cong_{\CC} (1357Q) \cong_{\CC}  (1357Q1)$) and $(1357QRS1)[t=0]$ ($\cong_{\CC} (2357D) \cong_{\CC} (2357D1) \cong_{\CC} (1357)[t=1]$) in an entirely analogous way.
\end{example}

\begin{example}\label{ex2}
In this example, we show how to find the nilsoliton metric for $(1357R)$ given by $(\RR^7,\mu)$ with
$$
\mu:=\left\{\begin{array}{l} [e_1, e_2] = e_3, [e_1, e_3] = e_5, [e_1, e_6] = e_7,\\{[e_2, e_3]} = e_6, [e_2, e_4] = e_6, [e_2, e_5] = e_7,[e_3, e_4] = e_7\end{array}\right.
$$
As we explained previously, $(1357R)$ must admit a such metric. Furthermore, there must exist a $g\in \mathrm{GL}_{7}(\RR)$ such that $\mm(g\cdot\mu)$ satisfies the Equation ($\ref{eq1}$).

It is easy to see that the pre-Einstein derivation of $(\RR^7,\mu)$ is equal to $\frac{4}{11}\diag(1,1,2,2,3,3,4)$ and as the pre-Einstein derivation of $g\cdot\mu$ must be positively proportional to the Einstein derivation, $\frac{4}{13}\diag(1,1,2,2,3,3,4)$, we can try to find a such $g$ in the group
$$
\mathrm{G}=\left\{g \in \mathrm{GL}_{7}(\RR) : g=\diag \left( \left(\begin{array}{ll} m_{1,1}&m_{1,2}\\ m_{2,1}&m_{2,2}\end{array}\right),\left(\begin{array}{ll} m_{3,3}&m_{3,4}\\ m_{4,3}&m_{4,4}\end{array}\right),\left(\begin{array}{ll} m_{5,5}&m_{5,6}\\ m_{6,5}&m_{6,6}\end{array}\right),m_{7,7} \right)\right\}
$$
which commutates with the Einstein derivation. By solving
\begin{eqnarray*}
\mm(g\cdot\mu) &=& \frac{1}{13}\diag(-7,-7,-3,-3,1,1,5)\\
\end{eqnarray*}
with $g\in\mathrm{G}$, we find
$$
g=\diag \left( \left(\begin{array}{cc} 1&\frac{\sqrt {13}}{2}\\ 1&-\frac{\sqrt {13}}{2}\end{array}\right),
\left(\begin{array}{cc}-\frac{2\sqrt {3}}{3}&-\frac{\sqrt {3}}{3}\\ 0&1\end{array}\right),
\left(\begin{array}{cc} -{\frac {2\sqrt {39}}{39}}&\frac{\sqrt {3}}{3}\\ -{\frac {2\sqrt {39}}{39}}&-\frac{\sqrt {3}}{3}\end{array}\right)
,-\frac{\sqrt {130}}{39} \right)
$$
which defines $\widetilde{\mu}=g\cdot\mu$
$$
\widetilde{\mu}:=\left\{\begin{array}{l}
[{e_1},{e_2}]={\frac {2\sqrt {39}}{39}}{e_3},[{e_1},{e_3}]={\frac {\sqrt {13}}{13}}{e_6},[{e_1},{e_4}]={\frac {\sqrt {39}}{39}}{e_5},[{e_1},{e_6}]={\frac {\sqrt {390}}{78}}{e_7},\\
{[{e_2},{e_3}]}={\frac {\sqrt {13}}{13}}{e_5},[{e_2},{e_4}]={\frac {\sqrt {39}}{39}}{e_6},[{e_2},{e_5}]=-{\frac {\sqrt {390}}{78}}{e_7},[{e_3},{e_4}]={\frac {\sqrt {390}}{78}}{e_7}
\end{array}\right.
$$
where the canonical inner product of $\RR^7$ defines a nilsoliton metric on $(\RR^7,\widetilde{\mu})$.
\end{example}

\begin{enumerate}\addtocounter{enumi}{75}
\item \textbf{(147A)}: $\dim(\Der)=15$, $\rank=3$, Betti Numbers $(4, 8, 9, 9, 8, 4, 1)$ ($\cong_{\RR} \ggo_{3.1(iii)}$).
  \newline Pre-Einstein der. $\phi=\frac{1}{2}\diag(1, 1, 1, 2, 2, 2, 3)$,
  \newline It does not admit

\item \textbf{(147A1)}: $\dim(\Der)=15$, $\rank=3$, Betti Numbers $(4, 8, 9, 9, 8, 4, 1)$ ($\cong_{\CC} \ggo_{3.1(iii)}$).
  \newline Pre-Einstein der. $\phi=\frac{1}{2}\diag(1, 1, 1, 2, 2, 2, 3)$,
  \newline It does not admit.

\item \textbf{(147B)}: $\dim(\Der)=12$, $\rank=2$, Betti Numbers $(4, 8, 10, 10, 8, 4, 1)$ ($\cong_{\RR} \ggo_{2.28}$).
  \newline Einstein der. $\phi=\frac{4}{35}\diag(4,6,5,10,9,8,14)$, $||\beta||^2=\frac{37}{35}\thickapprox1.057$
  \newline $[{e_1},{e_2}]={\frac {2\sqrt {35}}{35}}{e_4},[{e_1},{e_3}]={\frac {2\sqrt {35}}{35}}{e_5},[{e_1},{e_4}]={\frac {\sqrt {14}}{14}}{e_7},[{e_2},{e_6}]={\frac {\sqrt {14}}{14}}{e_7},[{e_3},{e_5}]={\frac {3\sqrt {70}}{70}}{e_7}$

\item \textbf{(147D)!}: $\dim(\Der)=15$, $\rank=2$, Betti Numbers $(3, 7, 9, 9, 7, 3, 1)$ ($\cong_{\RR} \ggo_{2.2}$!).
  \newline Pre-Einstein der. $\phi=\frac{1}{2}\diag(1,1,1,2,2,2,3)$,
  \newline It does not Einstein nilradical.

\item \textbf{(147E) $[0<t<1]$}: $\dim(\Der)=15$, $\rank=3$, Betti Numbers $(3, 7, 9, 9, 7, 3, 1)$ ($\cong_{\RR} \ggo_{3.1(i_{\lambda})}$).
  \newline Einstein der. $\phi=\frac{1}{2}\diag(1,1,1,2,2,2,3)$, $||\beta||^2=1$
  \newline $[{e_1},{e_2}]=\frac{\sqrt {2}}{4}\,\sqrt {1-t}{e_4},[{e_1},{e_3}]=-\frac{\sqrt {2}}{4}\,\sqrt {t}{e_6},[{e_1},{e_5}]=-\frac{\sqrt {2}}{4}{e_7},[{e_2},{e_3}]=\frac{\sqrt {2}}{4}{e_5},[{e_2},{e_6}]=\frac{\sqrt {2}}{4}\,\sqrt {t}{e_7},[{e_3},{e_4}]=\frac{\sqrt {2}}{4}\,\sqrt {1-t}{e_7}$

\item \textbf{(147E1) $[t>1]$}!: $\dim(\Der)=15$, $\rank=3$, Betti Numbers $(3, 7, 9, 9, 7, 3, 1)$ ($\cong_{\CC} \ggo_{3.1(i_{P(\lambda)})}$ $\cong_{\CC}$ (147E) $[\lambda]$ with $\lambda=\left(\frac{( 1-\sqrt{t^2-1}\sqrt{-1})}{t}\right)^2$ ).
  \newline Einstein der. $\phi=\frac{1}{2}\diag(1,1,1,2,2,2,3)$, $||\beta||^2=1$
  \newline It family admits nilsoliton metrics.


\item \textbf{(1457A)}: $\dim(\Der)=16$, $\rank=3$, Betti Numbers $(4, 6, 9, 9, 6, 4, 1)$ ($\cong_{\RR} \ggo_{3.17}$).
  \newline Einstein der. $\phi=\frac{5}{22}\diag(1,3,4,5,3,3,6)$, $||\beta||^2=\frac{21}{22}\thickapprox 0.9546$
  \newline $[{e_1},{e_2}]=\frac{\sqrt{66}}{22}{e_3},[{e_1},{e_3}]=\frac{\sqrt{77}}{22}{e_4},[{e_1},{e_4}]=\frac{\sqrt{33}}{22}{e_7},[{e_5},{e_6}]=\frac{\sqrt{66}}{22}{e_7}$

\item \textbf{(1457B)}: $\dim(\Der)=15$, $\rank=2$, Betti Numbers $(4, 6, 8, 8, 6, 4, 1)$ ($\cong_{\RR} \ggo_{2.30}$).
  \newline Einstein der. $\phi=\frac{4}{29}\diag(2,4,6,8,5,5,10)$, $||\beta||^2=\frac{27}{29}\thickapprox0.9310$
  \newline $[{e_1},{e_2}]={\frac {2\sqrt {29}}{29}}{e_3},[{e_1},{e_3}]={\frac {2\sqrt {29}}{29}}{e_4},[{e_1},{e_4}]={\frac {\sqrt {174}}{58}}{e_7},[{e_2},{e_3}]={\frac {\sqrt {174}}{58}}{e_7},[{e_5},{e_6}]={\frac {\sqrt {406}}{58}}{e_7}$


\item \textbf{(137A)}: $\dim(\Der)=14$, $\rank=3$, Betti Numbers $(4, 7, 8, 8, 7, 4, 1)$ ($\cong_{\RR} \ggo_{3.16}$).
  \newline Einstein der. $\phi=\frac{1}{14}\diag(5,8,5,8,13,13,18)$, $||\beta||^2=\frac{6}{7}\thickapprox 0.8571$
  \newline $[{e_1},{e_2}]=\frac{\sqrt{7}}{7}{e_5},[{e_1},{e_5}]=\frac{\sqrt{21}}{14}{e_7},[{e_3},{e_4}]=\frac{\sqrt{7}}{7}{e_6},[{e_3},{e_6}]=\frac{\sqrt{21}}{14}{e_7}$

\item \textbf{(137A1)}: $\dim(\Der)=14$, $\rank=3$, Betti Numbers $(4, 7, 8, 8, 7, 4, 1)$ ($\cong_{\CC} \ggo_{3.16}$).
  \newline Einstein der. $\phi=\frac{1}{14}\diag(5,5,8,8,13,13,18)$, $||\beta||^2=\frac{6}{7}\thickapprox0.8571$
  \newline $[{e_1},{e_3}]=\frac{\sqrt{14}}{14}{e_5},[{e_1},{e_4}]=\frac{\sqrt{14}}{14}{e_6},[{e_1},{e_5}]=\frac{\sqrt{21}}{14}{e_7},[{e_2},{e_3}]=-\frac{\sqrt{14}}{14}{e_6},[{e_2},{e_4}]=\frac{\sqrt{14}}{14}{e_5},[{e_2},{e_6}]=\frac{\sqrt{21}}{14}{e_7}$

\item \textbf{(137B)}: $\dim(\Der)=13$, $\rank=2$, Betti Numbers $(4, 7, 7, 7, 7, 4, 1)$ ($\cong_{\CC} \ggo_{2.23}$).
  \newline Pre-Einstein der. $\phi=\frac{4}{11}\diag(1, 2, 1, 2, 3, 3, 4)$,
  \newline It does not admit.

\item \textbf{(137B1)}: $\dim(\Der)=13$, $\rank=2$, Betti Numbers $(4, 7, 7, 7, 7, 4, 1)$ ($\cong_{\CC} \ggo_{2.23}$).
  \newline Pre-Einstein der. $\phi=\frac{4}{11}\diag(1, 1, 2, 2, 3, 3, 4)$,
  \newline It does not admit.

\item \textbf{(137C)}: $\dim(\Der)=15$, $\rank=3$, Betti Numbers $(4, 7, 8, 8, 7, 4, 1)$ ($\cong_{\RR} \ggo_{3.10}$).
  \newline Einstein der. $\phi=\frac{1}{22}\diag(7,12,11,16,19,23,30)$, $||\beta||^2=\frac{10}{11}\thickapprox0.9091$
  \newline $[{e_1},{e_2}]=\frac{\sqrt{11}}{11}{e_5},[{e_1},{e_4}]=\frac{\sqrt{11}}{11}{e_6},[{e_1},{e_6}]=\frac{\sqrt{55}}{22}{e_7},[{e_2},{e_3}]=\frac{\sqrt{11}}{11}{e_6},[{e_3},{e_5}]=-\frac{\sqrt{55}}{22}{e_7}$

\item \textbf{(137D)}: $\dim(\Der)=14$, $\rank=2$, Betti Numbers $(4, 7, 8, 8, 7, 4, 1)$ ($\cong_{\RR} \ggo_{2.1(iv)}$).
  \newline Pre-Einstein der. $\phi=\frac{2}{19}\diag(3, 6, 5, 8, 9, 11, 14)$,
  \newline It does not admit.


\item \textbf{(1357A)}: $\dim(\Der)=14$, $\rank=2$, Betti Numbers $(4, 7, 8, 8, 7, 4, 1)$ ($\cong_{\RR} \ggo_{2.1(iii)}$).
  \newline Einstein der. $\phi=\frac{2}{21}\diag(3,5,6,8,11,9,14)$, $||\beta||^2=\frac{19}{21}\thickapprox0.9048$
  \newline $[{e_1},{e_2}]=\frac{\sqrt{21}}{14}{e_4},[{e_1},{e_4}]=\frac{\sqrt{42}}{21}{e_5},[{e_1},{e_5}]=\frac{\sqrt{21}}{14}{e_7},[{e_2},{e_3}]=\frac{\sqrt{3}}{6}{e_5},
  [{e_2},{e_6}]=\frac{\sqrt{42}}{42}{e_7},[{e_3},{e_4}]=-\frac{\sqrt{3}}{6}{e_7}$

\item \textbf{(1357B)!}: $\dim(\Der)=14$, $\rank=2$, Betti Numbers $(4, 6, 7, 7, 6, 4, 1)$ ($\cong_{\RR} \ggo_{2.25}$!).
  \newline Einstein der. $\phi=\frac{5}{19}\diag(1,2,2,3,3,4,5)$, $||\beta||^2=\frac{17}{19}\thickapprox 0.8947$
  \newline $[{e_1},{e_3}]=\frac{\sqrt{38}}{19}{e_4},[{e_1},{e_4}]=-\frac{\sqrt{114}}{38}{e_6},[{e_1},{e_5}]=\frac{\sqrt{38}}{38}{e_6},[{e_1},{e_6}]=\frac{\sqrt{38}}{19}{e_7},[{e_2},{e_3}]=\frac{\sqrt{114}}{38}{e_6},[{e_2},{e_4}]=\frac{\sqrt{114}}{38}{e_7},[{e_2},{e_5}]=\frac{\sqrt{38}}{38}{e_7}$

\item \textbf{(1357C)!}: $\dim(\Der)=13$, $\rank=1$, Betti Numbers $(4, 6, 7, 7, 6, 4, 1)$ ($\cong_{\RR} \ggo_{1.3(v)}$!).
  \newline Pre-Einstein der. $\phi=\frac{5}{17}\diag(1, 2, 2, 3, 4, 3, 5)$,
  \newline It does not admit.

\item \textbf{(1357D)}: $\dim(\Der)=14$, $\rank=2$, Betti Numbers $(3, 6, 8, 8, 6, 3, 1)$ ($\cong_{\RR} \ggo_{2.1(ii)}$).
  \newline Einstein der. $\phi=\frac{2}{21}\diag(5,3,8,6,11,9,14)$, $||\beta||^2=\frac{19}{21}\thickapprox 0.9048$
  \newline $[{e_1},{e_2}]=\frac{\sqrt{21}}{14}{e_3},[{e_1},{e_6}]=\frac{\sqrt{21}}{14}{e_7},[{e_2},{e_3}]=\frac{\sqrt{42}}{21}{e_5},[{e_2},{e_4}]=-\frac{\sqrt{3}}{6}{e_6},[{e_2},{e_5}]=-\frac{\sqrt{42}}{42}{e_7},[{e_3},{e_4}]=-\frac{\sqrt{3}}{6}{e_7}$

\item \textbf{(1357E)}: $\dim(\Der)=14$, $\rank=2$, Betti Numbers $(3, 5, 8, 8, 5, 3, 1)$ ($\cong_{\RR} \ggo_{2.32}$).
  \newline Einstein der. $\phi=\frac{2}{29}\diag(10,3,13,8,16,11,19)$, $||\beta||^2=\frac{27}{29}\thickapprox0.9310$
  \newline $[{e_1},{e_2}]={\frac {\sqrt {406}}{58}}{e_3},[{e_2},{e_3}]={\frac {2\sqrt {29}}{29}}{e_5},[{e_2},{e_4}]={\frac {\sqrt {174}}{58}}{e_6},[{e_2},{e_5}]={\frac {\sqrt {174}}{58}}{e_7},[{e_4},{e_6}]={\frac {2\sqrt {29}}{29}}{e_7}$

\item \textbf{(1357F)}: $\dim(\Der)=13$, $\rank=1$, Betti Numbers $(3, 5, 7, 7, 5, 3, 1)$ ($\cong_{\RR} \ggo_{1.3(iii)}$).
  \newline Einstein der. $\phi=\frac{5}{19}\diag(2,1,3,2,4,3,5)$, $||\beta||^2=\frac{17}{19}\thickapprox 0.8947$
  \newline $[{e_1},{e_2}]={\frac {3\sqrt {76}}{76}}{e_3},[{e_1},{e_3}]={\frac {\sqrt {380}}{76}}{e_7},[{e_2},{e_3}]=\frac{\sqrt{38}}{19}{e_5},[{e_2},{e_4}]={\frac {\sqrt {380}}{76}}{e_6},[{e_2},{e_5}]=\frac{\sqrt{38}}{38}{e_7},[{e_4},{e_6}]=-{\frac {3\sqrt {76}}{76}}{e_7}$

\item \textbf{(1357F1)}: $\dim(\Der)=13$, $\rank=1$, Betti Numbers $(3, 5, 7, 7, 5, 3, 1)$ ($\cong_{\CC} \ggo_{1.3(iii)}$).
  \newline Einstein der. $\phi=\frac{5}{19}\diag(2,1,3,2,4,3,5)$, $||\beta||^2=\frac{17}{19}\thickapprox 0.8947$
  \newline $[{e_1},{e_2}]={\frac {3\sqrt {19}}{38}}{e_3},[{e_1},{e_3}]=\frac{\sqrt {95}}{38}{e_7},[{e_2},{e_3}]=\frac{\sqrt {38}}{19}{e_5},[{e_2},{e_4}]=\frac{\sqrt {95}}{38}{e_6},[{e_2},{e_5}]=\frac{\sqrt {38}}{38}{e_7},[{e_4},{e_6}]={\frac {3\sqrt {19}}{38}}{e_7}$

\item \textbf{(1357G)}: $\dim(\Der)=13$, $\rank=2$, Betti Numbers $(3, 6, 8, 8, 6, 3, 1)$ ($\cong_{\RR} \ggo_{2.31}$).
  \newline Einstein der. $\phi=\frac{1}{46}\diag(15,14,29,27,43,42,57)$, $||\beta||^2=\frac{39}{46}\thickapprox0.8478$
  \newline $[{e_1},{e_2}]=\frac{\sqrt{69}}{46}{e_3},[{e_1},{e_4}]=\frac{\sqrt{69}}{23}{e_6},[{e_1},{e_6}]={\frac {3\sqrt {23}}{46}}{e_7},[{e_2},{e_3}]=\frac{\sqrt{299}}{46}{e_5},[{e_2},{e_5}]={\frac {3\sqrt {23}}{46}}{e_7}$

\item \textbf{(1357H)!}: $\dim(\Der)=12$, $\rank=1$, Betti Numbers $(3, 6, 7, 7, 6, 3, 1)$ ($\cong_{\RR} \ggo_{1.2(iv)}$!).
  \newline Pre-Einstein der. $\phi=\frac{4}{11}\diag(1, 1, 2, 2, 3, 3, 4)$,
  \newline It does not admit.

\item \textbf{(1357I)}: $\dim(\Der)=12$, $\rank=2$, Betti Numbers $(3, 5, 7, 7, 5, 3, 1)$ ($\cong_{\RR} \ggo_{2.33}$).
  \newline Einstein der. $\phi=\frac{1}{41}\diag(18,10,28,15,38,33,48)$, $||\beta||^2=\frac{33}{41}\thickapprox0.8049$
  \newline $[{e_1},{e_2}]=\frac{\sqrt{123}}{41}{e_3},[{e_1},{e_4}]={\frac {3\sqrt {82}}{82}}{e_6},[{e_2},{e_3}]={\frac {\sqrt {902}}{82}}{e_5},[{e_2},{e_5}]=\frac{\sqrt{123}}{41}{e_7},[{e_4},{e_6}]={\frac {3\sqrt {82}}{82}}{e_7}$

\item \textbf{(1357J)}: $\dim(\Der)=11$, $\rank=1$, Betti Numbers $(3, 5, 6, 6, 5, 3, 1)$ ($\cong_{\RR} \ggo_{1.8}$).
  \newline Pre-Einstein der. $\phi=\frac{20}{139}\diag(4, 2, 6, 3, 8, 7, 10)$,
  \newline It does not admit.

\item \textbf{(1357L)!}: $\dim(\Der)=14$, $\rank=1$, Betti Numbers $(3, 6, 8, 8, 6, 3, 1)$ ($\cong_{\RR} \ggo_{1.3(ii)}$!).
  \newline Pre-Einstein der. $\phi=\frac{5}{17}\diag(1, 2, 3, 2, 4, 3, 5)$,
  \newline It does not admit.

\item\label{1357M} \textbf{(1357M)$[t\in\RR \setminus\{0,1\}]$}: $\dim(\Der)=14$, $\rank=2$, Betti Numbers $(3, 6, 8, 8, 6, 3, 1)$ ($\cong_{\RR} \ggo_{2.1(i_{t})}$).
  \newline Einstein der. $\phi=\frac{2}{21}\diag(3,5,8,6,11,9,14)$, $||\beta||^2=\frac{19}{21}\thickapprox0.9048$
  \newline It family admits nilsoliton metrics.

\item[(\ref{1357M}.1)] \textbf{(1357M)$[t=0]$}: $\dim(\Der)=14$, $\rank=2$,  Betti Numbers $(3, 6, 8, 8, 6, 3, 1)$ ($\cong_{\RR} (2357B)\cong_{\RR} \ggo_{2.1(i_{0})}$).
\newline Pre-Einstein der. $\phi=\frac{2}{19}\diag(3,5,8,6,11,9,14)$,
\newline It does not admit.

\item[(\ref{1357M}.2)] \textbf{(1357M)$[t=1]$}: $\dim(\Der)=14$, $\rank=2$, Betti Numbers $(3, 6, 8, 8, 6, 3, 1)$ ($\cong_{\RR} \ggo_{2.1(i_{1})}$).
  \newline Einstein der. $\phi=\frac{2}{21}\diag(3,5,8,6,11,9,14)$, $||\beta||^2=\frac{19}{21}\thickapprox0.9048$
  \newline $[{e_1},{e_2}]=\frac{\sqrt{42}}{42}{e_3},[{e_1},{e_3}]=\frac{\sqrt{42}}{21}{e_5},[{e_1},{e_4}]=\frac{\sqrt{3}}{6}{e_6},[{e_1},{e_5}]=\frac{\sqrt{21}}{14}{e_7},[{e_2},{e_4}]=\frac{\sqrt{3}}{6}{e_5},[{e_2},{e_6}]=\frac{\sqrt{21}}{14}{e_7}$

 \item\label{1357N} \textbf{(1357N)$[t\in\RR\setminus\{0\}]$}!: $\dim(\Der)=13$, $\rank=1$, Betti Numbers $(3, 5, 7, 7, 5, 3, 1)$ ($\cong_{\RR} \ggo_{1.3(i_{\lambda})}$!).
  \newline Einstein der. $\phi=\frac{5}{19}\diag(1,2,3,2,4,3,5)$, $||\beta||^2=\frac{17}{19}\thickapprox0.8947$
  \newline It family admits nilsoliton metrics.

 \item[(\ref{1357N}.1)] \textbf{(1357N)$[t=0]$}!: $\dim(\Der)=13$, $\rank=1$, Betti Numbers $(3, 5, 7, 7, 5, 3, 1)$ ($\cong_{\RR} \ggo_{1.3(i_{0})}$!).
  \newline Pre-Einstein der. $\phi=\frac{5}{17}\diag(1,2,3,2,4,3,5)$,
  \newline It does not admit.

\item \textbf{(1357O)}: $\dim(\Der)=13$, $\rank=2$, Betti Numbers $(3, 6, 8, 8, 6, 3, 1)$ ($\cong_{\RR} \ggo_{2.41}$).
  \newline Einstein der. $\phi=\frac{1}{8}\diag(3,2,5,6,8,7,10)$, $||\beta||^2=\frac{7}{8}\thickapprox0.8750$
  \newline $[{e_1},{e_2}]=\frac{1}{4}{e_3},[{e_1},{e_3}]=\frac{\sqrt{6}}{8}{e_5},[{e_1},{e_6}]=\frac{\sqrt{6}}{8}{e_7},[{e_2},{e_3}]=\frac{\sqrt{6}}{8}{e_6},[{e_2},{e_4}]=\frac{1}{4}{e_5},[{e_2},{e_5}]=\frac{\sqrt{6}}{8}{e_7}
$

\item \textbf{(1357P)}: $\dim(\Der)=12$, $\rank=1$, Betti Numbers $(3, 6, 7, 7, 6, 3, 1)$ ($\cong_{\RR} \ggo_{1.2(i_{0})}$).
  \newline Einstein der. $\phi=\frac{4}{13}\diag(1,1,2,2,3,3,4)$, $||\beta||^2=\frac{11}{13}\thickapprox 0.8462$
  \newline $[{e_1},{e_2}]=\frac{\sqrt{65}}{26}{e_3},[{e_1},{e_3}]=\frac{\sqrt{13}}{13}{e_5},[{e_1},{e_5}]=\frac{\sqrt{65}}{26}{e_7},[{e_2},{e_3}]=\frac{\sqrt{13}}{13}{e_6},[{e_2},{e_4}]=\frac{\sqrt{39}}{26}{e_5},[{e_2},{e_6}]=\frac{\sqrt{26}}{26}{e_7},[{e_3},{e_4}]=\frac{\sqrt{39}}{26}{e_7}$

\item \textbf{(1357P1)}: $\dim(\Der)=12$, $\rank=1$, Betti Numbers $(3, 6, 7, 7, 6, 3, 1)$ ($\cong_{\CC} \ggo_{1.2(i_{0})}$).
  \newline Einstein der. $\phi=\frac{4}{13}\diag(1,1,2,2,3,3,4)$, $||\beta||^2=\frac{11}{13}\thickapprox 0.8462$
  \newline $[{e_1},{e_2}]=\frac{\sqrt{65}}{26}{e_3},[{e_1},{e_3}]=\frac{\sqrt{13}}{13}{e_5},[{e_1},{e_5}]=\frac{\sqrt{65}}{26}{e_7},[{e_2},{e_3}]=\frac{\sqrt{13}}{13}{e_6},[{e_2},{e_4}]=\frac{\sqrt{39}}{26}{e_5},[{e_2},{e_6}]=-\frac{\sqrt{26}}{26}{e_7},[{e_3},{e_4}]=\frac{\sqrt{39}}{26}{e_7}$

\item \textbf{(1357Q)!}: $\dim(\Der)=12$, $\rank=1$, Betti Numbers $(3, 6, 8, 8, 6, 3, 1)$ ($\cong_{\RR} \ggo_{1.2(i_{1})}$).
  \newline Einstein der. $\phi=\frac{4}{13}\diag(1,1,2,2,3,3,4)$, $||\beta||^2=\frac{11}{13}\thickapprox0.8462$
  \newline $[{e_1},{e_2}]=\frac{\sqrt{26}}{26}{e_3},[{e_1},{e_3}]=\frac{\sqrt{13}}{13}{e_5},[{e_1},{e_4}]=\frac{\sqrt{39}}{26}{e_6},[{e_1},{e_5}]=\frac{\sqrt{65}}{26}{e_7},[{e_2},{e_3}]=\frac{\sqrt{13}}{13}{e_6},[{e_2},{e_4}]=\frac{\sqrt{39}}{26}{e_5},[{e_2},{e_6}]=\frac{\sqrt{65}}{26}{e_7}
$

\item \textbf{(1357Q1)!}: $\dim(\Der)=12$, $\rank=1$, Betti Numbers $(3, 6, 8, 8, 6, 3, 1)$ ($\cong_{\CC} \ggo_{1.2(i_{1})}$).
  \newline Einstein der. $\phi=\frac{4}{13}\diag(1,1,2,2,3,3,4)$, $||\beta||^2=\frac{11}{13}\thickapprox 0.8462$
  \newline $[{e_1},{e_2}]=\frac{\sqrt{26}}{26}{e_3},[{e_1},{e_3}]=\frac{\sqrt{13}}{13}{e_5},[{e_1},{e_4}]=\frac{\sqrt{39}}{26}{e_6},[{e_1},{e_6}]=\frac{\sqrt{65}}{26}{e_7},[{e_2},{e_3}]=\frac{\sqrt{13}}{13}{e_6},[{e_2},{e_4}]=\frac{\sqrt{39}}{26}{e_5},[{e_2},{e_5}]=\frac{\sqrt{65}}{26}{e_7}$

\item \textbf{(1357R)!}: $\dim(\Der)=13$, $\rank=2$, Betti Numbers $(3, 6, 7, 7, 6, 3, 1)$ ($\cong_{\RR} \ggo_{2.37}$!).
  \newline Einstein der. $\phi=\frac{4}{11}\diag(1, 1, 2, 2, 3, 3, 4)$,
  $||\beta||^2={\frac {11}{13}}\thickapprox0.8462$
  \newline $ [{e_1},{e_2}]={\frac {2\sqrt {39}}{39}}{e_3},[{e_1},{e_3}]={\frac {\sqrt {13}}{13}}{e_6},[{e_1},{e_4}]={\frac {\sqrt {39}}{39}}{e_5},[{e_1},{e_6}]={\frac {\sqrt {390}}{78}}{e_7},[{e_2},{e_3}]={\frac {\sqrt {13}}{13}}{e_5},[{e_2},{e_4}]={\frac {\sqrt {39}}{39}}{e_6},[{e_2},{e_5}]=-{\frac {\sqrt {390}}{78}}{e_7},[{e_3},{e_4}]={\frac {\sqrt {390}}{78}}{e_7}$

\item\label{1357QRS1} \textbf{(1357QRS1) $[t\in\RR\setminus\{-1,0,1\}]$}: $\dim(\Der)=12$, $\rank=1$, Betti Numbers $(3, 6, 7, 7, 6, 3, 1)$ ($\cong_{\CC} \ggo_{1.2(i_{\lambda})}$ with $\lambda=\frac{1}{t}$).
  \newline Einstein der. $\phi=\frac{4}{13}\diag(1,1,2,2,3,3,4)$, $||\beta||^2={\frac {11}{13}}\thickapprox  0.8462$
  \newline It family admits nilsoliton metrics.

\item[(\ref{1357QRS1}.1)] \textbf{(1357QRS1) $[t=-1]$}: $\dim(\Der)=13$, $\rank=2$, Betti Numbers $(3, 6, 7, 7, 6, 3, 1)$ ($\cong_{\CC} \ggo_{2.37} \cong_{\CC}$ (1357R))
  \newline Einstein der. $\phi=\frac{4}{13}\diag(1,1,2,2,3,3,4)$, $||\beta||^2={\frac {11}{13}}\thickapprox0.8462$
  \newline $[{e_1},{e_2}]={\frac {2\sqrt {39}}{39}}{e_3},[{e_1},{e_3}]=\frac{\sqrt{13}}{13}{e_5},[{e_1},{e_4}]=\frac{\sqrt{39}}{39}{e_6},[{e_1},{e_5}]={\frac {\sqrt {390}}{78}}{e_7},[{e_2},{e_3}]=-\frac{\sqrt{13}}{13}{e_6},[{e_2},{e_4}]=\frac{\sqrt{39}}{39}{e_5},[{e_2},{e_6}]=-{\frac {\sqrt {390}}{78}}{e_7},[{e_3},{e_4}]={\frac {\sqrt {390}}{78}}{e_7}$

\item[(\ref{1357QRS1}.2)] \textbf{(1357QRS1)$[t=1]$}: $\dim(\Der)=12$, $\rank=1$, Betti Numbers $(3, 6, 8, 8, 6, 3, 1)$ ($\cong_{\CC} \ggo_{1.2(i_{1})} \cong_{\CC}$ (1357Q)).
  \newline Einstein der. $\phi=\frac{4}{13}\diag(1,1,2,2,3,3,4)$, $||\beta||^2=\frac{11}{13}\thickapprox0.8462$
  \newline $[{e_1},{e_2}]=\frac{\sqrt{26}}{26}{e_3},[{e_1},{e_3}]=\frac{\sqrt{13}}{13}{e_5},[{e_1},{e_4}]=\frac{\sqrt{39}}{26}{e_6},[{e_1},{e_5}]=\frac{\sqrt{65}}{26}{e_7},[{e_2},{e_3}]=-\frac{\sqrt{13}}{13}{e_6},[{e_2},{e_4}]=\frac{\sqrt{39}}{26}{e_5},[{e_2},{e_6}]=\frac{\sqrt{65}}{26}{e_7}$

\item[(\ref{1357QRS1}.3)] \textbf{(1357QRS1) $[t=0]$}: $\dim(\Der)=12$, $\rank=1$, Betti Numbers $(3, 6, 7, 7, 6, 3, 1)$ ($\cong_{\CC} \ggo_{1.2(iii)} \cong_{\CC}$ (2357D)).
  \newline Pre-Einstein der. $\phi=\frac{4}{11}\diag(1,1,2,2,3,3,4)$,
  \newline It does not admit.

\item\label{1357S} \textbf{(1357S)! $[t\in\RR\setminus\{0,1\}]$}: $\dim(\Der)=12$, $\rank=1$, Betti Numbers $(3, 6, 7, 7, 6, 3, 1)$ ($\cong_{\CC} \ggo_{1.2(i_{P(\lambda)})} \cong_{\CC}$ (1357QRS1) $[u]$ with $u=\frac{2\sqrt{t}+t+1}{t-1}$).
  \newline Einstein der. $\phi=\frac{4}{13}\diag(1,1,2,2,3,3,4)$, $||\beta||^2={\frac {11}{13}}\thickapprox 0.8462$
  \newline It family admits nilsoliton metrics.

\item[(\ref{1357S}.1)] \textbf{(1357S)! $[t=0]$}: $\dim(\Der)=12$, $\rank=1$, Betti Numbers $(3, 6, 7, 7, 6, 3, 1)$ ($\cong_{\RR} \ggo_{1.2(ii)}$!).
  \newline Pre-Einstein der. $\phi=\frac{4}{11}\diag(1,1,2,2,3,3,4)$,
  \newline It does not admit.

\item[(\ref{1357S}.2)] \textbf{(1357S)! $[t=1]$}: $\dim(\Der)=12$, $\rank=1$, Betti Numbers $(3, 6, 7, 7, 6, 3, 1)$ ($\cong_{\CC} \ggo_{1.2(iii)} \cong_{\CC}$ (2357D)).
  \newline Pre-Einstein der. $\phi=\frac{4}{11}\diag(1,1,2,2,3,3,4)$,
  \newline It does not admit.


\item \textbf{(13457A)}: $\dim(\Der)=14$, $\rank=2$, Betti Numbers $(3, 5, 7, 7, 5, 3, 1)$ ($\cong_{\RR} \ggo_{2.16}$).
  \newline Einstein der. $\phi=\frac{1}{29}\diag(5,17,22,27,32,20,37)$, $||\beta||^2=\frac{27}{29}\thickapprox 0.9310$
  \newline $[{e_1},{e_2}]={\frac {\sqrt {174}}{58}}{e_3},[{e_1},{e_3}]={\frac {2\sqrt {29}}{29}}{e_4},[{e_1},{e_4}]={\frac {2\sqrt {29}}{29}}{e_5},[{e_1},{e_5}]={\frac {\sqrt {174}}{58}}{e_7},[{e_2},{e_6}]={\frac {\sqrt {406}}{58}}{e_7}$

\item \textbf{(13457B)!}: $\dim(\Der)=13$, $\rank=1$, Betti Numbers $(3, 5, 7, 7, 5, 3, 1)$ ($\cong_{\RR} \ggo_{1.15}$!).
  \newline Einstein der. $\phi=\frac{15}{82}\diag(1,3,4,4,5,6,7)$, $||\beta||^2=\frac{38}{41}\thickapprox 0.9268$
  \newline $[{e_1},{e_2}]={\frac {2\sqrt {4305}}{861}}{e_3}+{\frac {\sqrt {111930}}{1722}}{e_4},[{e_1},{e_3}]={\frac {2\sqrt {22386}}{861}}{e_5},[{e_1},{e_4}]={\frac {5\sqrt {861}}{1722}}{e_5},[{e_1},{e_5}]={\frac {\sqrt {902}}{82}}{e_6},[{e_1},{e_6}]=\frac{\sqrt{82}}{41}{e_7},[{e_2},{e_4}]={\frac {\sqrt {861}}{82}}{e_7}$

\item \textbf{(13457C)}: $\dim(\Der)=12$, $\rank=2$, Betti Numbers $(3, 4, 4, 4, 4, 3, 1)$ ($\cong_{\RR} \ggo_{2.10}$).
  \newline Pre-Einstein der. $\phi=\frac{1}{5}\diag(1, 2, 3, 4, 5, 6, 7)$,
  \newline It does not admit.

\item \textbf{(13457D)!}: $\dim(\Der)=12$, $\rank=1$, Betti Numbers $(3, 5, 7, 7, 5, 3, 1)$ ($\cong_{\RR} \ggo_{1.12}$!).
  \newline Einstein der. $\phi=\frac{25}{124}\diag(1,2,4,3,4,5,6)$, $||\beta||^2=\frac{107}{124}\thickapprox  0.8629$
  \newline $[{e_1},{e_2}]={\frac {\sqrt {930}}{124}}{e_4},[{e_1},{e_3}]=\frac{\sqrt{31}}{31}{e_6},[{e_1},{e_4}]={\frac {\sqrt {341}}{62}}{e_3}+{\frac {\sqrt {62}}{62}}{e_5},[{e_1},{e_5}]={\frac {\sqrt {682}}{124}}{e_6},[{e_1},{e_6}]={\frac {\sqrt {341}}{62}}{e_7},[{e_2},{e_3}]={\frac {\sqrt {1302}}{124}}{e_7},[{e_2},{e_4}]={\frac {\sqrt {1302}}{124}}{e_6}
$

\item \textbf{(13457E)}: $\dim(\Der)=11$, $\rank=1$, Betti Numbers $(3, 4, 4, 4, 4, 3, 1)$ ($\cong_{\RR} \ggo_{1.1(vi)}$).
  \newline Pre-Einstein der. $\phi=\frac{1}{5}\diag(1, 2, 3, 4, 5, 6, 7)$,
  \newline It does not admit.

\item \textbf{(13457F)}: $\dim(\Der)=11$, $\rank=1$, Betti Numbers $(2, 4, 7, 7, 4, 2, 1)$ ($\cong_{\RR} \ggo_{1.10}$).
  \newline Einstein der. $\phi=\frac{45}{446}\diag(2,3,5,7,9,8,11)$, $||\beta||^2={\frac {353}{446}}\thickapprox  0.7915$
  \newline $[{e_1},{e_2}]={\frac {\sqrt {31220}}{892}}{e_3},[{e_1},{e_3}]={\frac {\sqrt {15164}}{446}}{e_4},[{e_1},{e_4}]={\frac {\sqrt {23638}}{446}}{e_5},[{e_1},{e_5}]={\frac {3\sqrt {1338}}{446}}{e_7},[{e_2},{e_3}]={\frac {\sqrt {84740}}{892}}{e_6},[{e_2},{e_6}]={\frac {\sqrt {4906}}{223}}{e_7}$

\item \textbf{(13457G)!!}: $\dim(\Der)=11$, $\rank=1$, Betti Numbers $(2, 3, 4, 4, 3, 2, 1)$ ($\cong_{\RR} \ggo_{1.03}$!!).
  \newline Pre-Einstein der. $\phi=\frac{2}{3}\diag(0, 1, 1, 1, 1, 2, 2)$,
  \newline It does not admit; $\phi \ngtr0$

\item \textbf{(13457I)!!}: $\dim(\Der)=10$, $\rank=0$, Betti Numbers $(2, 3, 4, 4, 3, 2, 1)$ ($\cong_{\RR} \ggo_{0.7}$!!).
  \newline Pre-Einstein der. $\phi=0$,
  \newline It does not admit; $\phi \ngtr0$


\item \textbf{(12457A)}: $\dim(\Der)=13$, $\rank=2$, Betti Numbers $(3, 5, 7, 7, 5, 3, 1)$ ($\cong_{\RR} \ggo_{2.15}$).
  \newline Pre-Einstein der. $\phi=\frac{1}{6}\diag(1,3,4,5,3,6,7)$, $||\beta||^2=\frac{5}{6}\thickapprox  0.8333$
  \newline $[{e_1},{e_2}]=\frac{\sqrt{3}}{6}{e_3},[{e_1},{e_3}]=\frac{\sqrt{3}}{6}{e_4},[{e_1},{e_4}]=\frac{\sqrt{3}}{6}{e_6},[{e_1},{e_6}]=\frac{\sqrt{3}}{6}{e_7},[{e_2},{e_5}]=\frac{\sqrt{3}}{6}{e_6},[{e_3},{e_5}]=\frac{\sqrt{3}}{6}{e_7}$

\item \textbf{(12457B)!!}: $\dim(\Der)=12$, $\rank=1$, Betti Numbers $(3, 5, 7, 7, 5, 3, 1)$ ($\cong_{\RR} \ggo_{1.01(ii)}$!!).
  \newline Pre-Einstein der. $\phi=\diag(0, 1, 1, 1, 0, 1, 1)$,
  \newline It does not admit; $\phi \ngtr0$

\item \textbf{(12457C)}: $\dim(\Der)=12$, $\rank=2$, Betti Numbers $(3, 4, 4, 4, 4, 3, 1)$ ($\cong_{\RR} \ggo_{2.13}$).
  \newline Einstein der. $\phi=\frac{1}{86}\diag(16,21,37,53,48,69,90)$, $||\beta||^2=\frac{30}{43}\thickapprox 0.6977$
  \newline $[{e_1},{e_2}]=\frac{\sqrt{129}}{43}{e_3},[{e_1},{e_3}]=\frac{\sqrt{215}}{43}{e_4},[{e_1},{e_4}]=\frac{\sqrt{129}}{43}{e_6},[{e_2},{e_5}]=\frac{\sqrt{129}}{43}{e_6},[{e_2},{e_6}]={\frac {\sqrt {645}}{86}}{e_7},[{e_3},{e_4}]=-{\frac {\sqrt {645}}{86}}{e_7}$

\item \textbf{(12457D)}: $\dim(\Der)=11$, $\rank=1$, Betti Numbers $(3, 4, 4, 4, 4, 3, 1)$ ($\cong_{\RR} \ggo_{1.20}$).
  \newline Pre-Einstein der. $\phi=\frac{8}{47}\diag(2, 1, 3, 5, 6, 7, 8)$,
  \newline It does not admit.

\item \textbf{(12457E)!}: $\dim(\Der)=11$, $\rank=1$, Betti Numbers $(3, 5, 6, 6, 5, 3, 1)$ ($\cong_{\RR} \ggo_{1.11}$!).
  \newline Einstein der. $\phi=\frac{6}{31}\diag(1,2,3,3,4,5,6)$, $||\beta||^2=\frac{25}{31}\thickapprox 0.8064$
  \newline $[{e_1},{e_2}]={\frac {7\sqrt {1767}}{1767}}{e_3}+{\frac {\sqrt {465}}{93}}{e_4},[{e_1},{e_3}]={\frac {\sqrt {651}}{93}}{e_5},[{e_1},{e_4}]={\frac {\sqrt {61845}}{1767}}{e_5},[{e_1},{e_5}]=\frac{\sqrt {62}}{31}{e_6},[{e_1},{e_6}]={\frac {\sqrt {90706}}{1178}}{e_7},[{e_2},{e_4}]={\frac {4\sqrt {1767}}{589}}{e_6},[{e_2},{e_5}]={\frac {\sqrt {64790}}{1178}}{e_7},[{e_3},{e_4}]={\frac {\sqrt {90706}}{1178}}{e_7}$

\item \textbf{(12457F)!}: $\dim(\Der)=11$, $\rank=1$, Betti Numbers $(3, 4, 4, 4, 4, 3, 1)$ ($\cong_{\RR} \ggo_{1.21}$!).
  \newline Pre-Einstein der. $\phi=\frac{25}{113}\diag(1, 2, 3, 4, 3, 5, 7)$,
  \newline It does not admit.

\item \textbf{(12457G)!!}: $\dim(\Der)=10$, $\rank=0$, Betti Numbers $(3, 4, 4, 4, 4, 3, 1)$ ($\cong_{\RR} \ggo_{0.8}$!!).
  \newline Pre-Einstein der. $\phi=0$,
  \newline It does not admit; $\phi \ngtr0$.

\item \textbf{(12457H)}: $\dim(\Der)=12$, $\rank=2$, Betti Numbers $(2, 3, 6, 6, 3, 2, 1)$ ($\cong_{\RR} \ggo_{2.5}$).
  \newline Einstein der. $\phi=\frac{1}{7}\diag(1,2,3,4,5,6,7)$, $||\beta||^2=\frac{5}{7}\thickapprox 0.7143$
  \newline $[{e_1},{e_2}]=\frac{\sqrt{14}}{14}{e_3},[{e_1},{e_3}]=\frac{\sqrt{14}}{14}{e_4},[{e_1},{e_5}]=\frac{\sqrt{14}}{14}{e_6},[{e_1},{e_6}]=\frac{\sqrt{14}}{14}{e_7},[{e_2},{e_3}]=\frac{\sqrt{14}}{14}{e_5},[{e_2},{e_4}]=\frac{\sqrt{14}}{14}{e_6},[{e_3},{e_4}]=\frac{\sqrt{14}}{14}{e_7}$

\item \textbf{(12457I)}: $\dim(\Der)=11$, $\rank=1$, Betti Numbers $(2, 3, 6, 6, 3, 2, 1)$ ($\cong_{\RR} \ggo_{1.1(iv)}$).
  \newline Pre-Einstein der. $\phi=\frac{1}{5}\diag(1, 2, 3, 4, 5, 6, 7)$,
  \newline It does not admit.

\item \textbf{(12457J)!!}: $\dim(\Der)=10$, $\rank=0$, Betti Numbers $(2, 3, 5, 5, 3, 2, 1)$ ($\cong_{\RR} \ggo_{0.6}$!!).
  \newline Pre-Einstein der. $\phi=0$,
  \newline It does not admit; $\phi \ngtr0$.

\item \textbf{(12457J1)!!}: $\dim(\Der)=10$, $\rank=0$, Betti Numbers $(2, 3, 5, 5, 3, 2, 1)$ ($\cong_{\CC} \ggo_{0.6}$!!).
  \newline Pre-Einstein der. $\phi=0$,
  \newline It does not admit; $\phi \ngtr0$.

\item \textbf{(12457K)!!}: $\dim(\Der)=11$, $\rank=1$, Betti Numbers $(2, 3, 5, 5, 3, 2, 1)$ ($\cong_{\RR} \ggo_{1.02}$!!).
  \newline Pre-Einstein der. $\phi=\frac{1}{2}\diag(1, 0, 1, 2, 1, 2, 3)$,
  \newline It does not admit; $\phi \ngtr0$.

\item \textbf{(12457L)!}: $\dim(\Der)=11$, $\rank=1$, Betti Numbers $(2, 3, 4, 4, 3, 2, 1)$ ($\cong_{\CC} \ggo_{1.17}$!).
  \newline Einstein der. $\phi=\frac{19}{94}\diag(1,1,2,3,3,4,5)$, $||\beta||^2=\frac{65}{94}\thickapprox  0.6915$
  \newline $[{e_1},{e_2}]={\frac {\sqrt {611}}{94}}{e_3},[{e_1},{e_3}]={\frac {\sqrt {235}}{94}}{e_4},
  [{e_1},{e_4}]=-{\frac {\sqrt {611}}{94}}{e_6},[{e_1},{e_6}]={\frac {\sqrt {705}}{94}}{e_7},[{e_2},{e_3}]=\frac{\sqrt{235}}{47}{e_5},
  [{e_2},{e_5}]={\frac {\sqrt {611}}{94}}{e_6},[{e_3},{e_5}]={\frac {\sqrt {705}}{94}}{e_7}$

\item \textbf{(12457L1)}: $\dim(\Der)=11$, $\rank=1$, Betti Numbers $(2, 3, 4, 4, 3, 2, 1)$ ($\cong_{\RR} \ggo_{1.17}$!).
  \newline Einstein der. $\phi=\frac{19}{94}\diag(1,1,2,3,3,4,5)$, $||\beta||^2=\frac{65}{94}\thickapprox  0.6915$
  \newline $[{e_1},{e_2}]={\frac {\sqrt {611}}{94}}{e_3},[{e_1},{e_3}]={\frac {\sqrt {235}}{94}}{e_4},
  [{e_1},{e_4}]=-{\frac {\sqrt {611}}{94}}{e_6},[{e_1},{e_6}]={\frac {\sqrt {705}}{94}}{e_7},[{e_2},{e_3}]={\frac {\sqrt {235}}{47}}{e_5},
  [{e_2},{e_5}]=-{\frac {\sqrt {611}}{94}}{e_6},[{e_3},{e_5}]=-{\frac {\sqrt {705}}{94}}{e_7}$

\item \textbf{(12457N)[t]!!}: $\dim(\Der)=10$, $\rank=0$, Betti Numbers $(2, 3, 4, 4, 3, 2, 1)$ ($\cong_{\CC} \ggo_{0.4_{P(\lambda)}}$!!).
  \newline Pre-Einstein der. $\phi=0$,
  \newline It does not admit; $\phi \ngtr0$.

\item \textbf{(12457N1)!!}: $\dim(\Der)=10$, $\rank=0$, Betti Numbers $(2, 3, 4, 4, 3, 2, 1)$ ($\cong_{\CC} \ggo_{0.4_{\lambda_{0}}}$!!).
  \newline Einstein der. $\phi=0$,
  \newline It does not admit; $\phi \ngtr0$.

\item \textbf{(12457N2)[t $\geq 0$]!!}: $\dim(\Der)=10$, $\rank=0$, Betti Numbers $(2, 3, 4, 4, 3, 2, 1)$ ($\cong_{\CC} \ggo_{0.4(\lambda)}$!!).
  \newline Pre-Einstein der. $\phi=0$,
  \newline It does not admit; $\phi \ngtr0$


\item \textbf{(12357A)}: $\dim(\Der)=12$, $\rank=2$, Betti Numbers $(3, 4, 4, 4, 4, 3, 1)$ ($\cong_{\RR} \ggo_{2.14}$).
  \newline Einstein der. $\phi=\frac{1}{7}\diag(1,3,2,4,5,6,7)$, $||\beta||^2=\frac{5}{7}\thickapprox 0.7143$
  \newline $[{e_1},{e_2}]=\frac{\sqrt{14}}{14}{e_4},[{e_1},{e_4}]=\frac{\sqrt{14}}{14}{e_5},[{e_1},{e_5}]=\frac{\sqrt{14}}{14}{e_6},[{e_1},{e_6}]=\frac{\sqrt{14}}{14}{e_7},[{e_2},{e_3}]=\frac{\sqrt{14}}{14}{e_5},[{e_3},{e_4}]=-\frac{\sqrt{14}}{14}{e_6},[{e_3},{e_5}]=-\frac{\sqrt{14}}{14}{e_7}$

\item \textbf{(12357B)!!}: $\dim(\Der)=11$, $\rank=1$, Betti Numbers $(3, 4, 4, 4, 4, 3, 1)$ ($\cong_{\RR} \ggo_{1.01(i)}$!!).
  \newline Pre-Einstein der. $\phi=\diag(0, 1, 0, 1, 1, 1, 1)$,
  \newline It does not admit; $\phi \ngtr0$

\item \textbf{(12357B1)!!}: $\dim(\Der)=11$, $\rank=1$, Betti Numbers $(3, 4, 4, 4, 4, 3, 1)$ ($\cong_{\CC} \ggo_{1.01(i)}$!!).
  \newline Pre-Einstein der. $\phi=\diag(0, 1, 0, 1, 1, 1, 1)$,
  \newline  It does not admit; $\phi \ngtr0$

\item \textbf{(12357C)}: $\dim(\Der)=10$, $\rank=1$, Betti Numbers $(3, 4, 4, 4, 4, 3, 1)$ ($\cong_{\RR} \ggo_{1.1(v)}$).
  \newline Pre-Einstein der. $\phi=\frac{1}{5}\diag(1, 3, 2, 4, 5, 6, 7)$,
  \newline It does not admit.


\item \textbf{(123457A)}: $\dim(\Der)=13$, $\rank=2$, Betti Numbers $(2, 4, 6, 6, 4, 2, 1)$ ($\cong_{\RR} \ggo_{2.3}$).
  \newline Einstein der. $\phi=\frac{2}{35}\diag(1,16,17,18,19,20,21)$, $||\beta||^2=\frac{37}{35}\thickapprox1.057$
  \newline $[{e_1},{e_2}]=\frac{\sqrt{14}}{14}{e_3},[{e_1},{e_3}]={\frac {2\sqrt {35}}{35}}{e_4},[{e_1},{e_4}]={\frac {3\sqrt {70}}{70}}{e_5},[{e_1},{e_5}]={\frac {2\sqrt {35}}{35}}{e_6},[{e_1},{e_6}]=\frac{\sqrt{14}}{14}{e_7}$

\item \textbf{(123457B)}: $\dim(\Der)=12$, $\rank=1$, Betti Numbers $(2, 4, 6, 6, 4, 2, 1)$ ($\cong_{\RR} \ggo_{1.6}$).
  \newline Einstein der. $\phi=\frac{5}{38}\diag(1,4,5,6,7,8,9)$, $||\beta||^2=\frac{17}{19}\thickapprox0.8947$
  \newline $[{e_1},{e_2}]={\frac {\sqrt {380}}{76}}{e_3},[{e_1},{e_3}]={\frac {\sqrt {380}}{76}}{e_4},[{e_1},{e_4}]={\frac {3\sqrt {76}}{76}}{e_5},[{e_1},{e_5}]=\frac{\sqrt{38}}{19}{e_6},[{e_1},{e_6}]=\frac{\sqrt{38}}{38}{e_7},[{e_2},{e_3}]={\frac {3\sqrt {76}}{76}}{e_7}$

\item \textbf{(123457C)}: $\dim(\Der)=11$, $\rank=1$, Betti Numbers $(2, 3, 4, 4, 3, 2, 1)$ ($\cong_{\RR} \ggo_{1.1(ii)}$).
  \newline Pre-Einstein der. $\phi=\frac{1}{5}\diag(1, 2, 3, 4, 5, 6, 7)$,
  \newline It does not admit.

\item \textbf{(123457D)}: $\dim(\Der)=12$, $\rank=1$, Betti Numbers $(2, 4, 6, 6, 4, 2, 1)$ ($\cong_{\RR} \ggo_{1.4}$).
  \newline Einstein der. $\phi=\frac{17}{122}\diag(1,3,4,5,6,7,8)$, $||\beta||^2=\frac{50}{61}\thickapprox 0.8197$
  \newline $[{e_1},{e_2}]={\frac {\sqrt {610}}{122}}{e_3},[{e_1},{e_3}]={\frac {\sqrt {1281}}{122}}{e_4},[{e_1},{e_4}]={\frac {3\sqrt {122}}{122}}{e_5},[{e_1},{e_5}]={\frac {\sqrt {244}}{61}}{e_6},[{e_1},{e_6}]={\frac {3\sqrt {122}}{122}}{e_7},[{e_2},{e_3}]={\frac {\sqrt {5124}}{244}}{e_6},[{e_2},{e_4}]={\frac {3\sqrt {122}}{122}}{e_7}$

\item \textbf{(123457E)!!}: $\dim(\Der)=11$, $\rank=0$, Betti Numbers $(2, 4, 6, 6, 4, 2, 1)$ ($\cong_{\RR} \ggo_{0.3}$!!).
  \newline Pre-Einstein der. $\phi=0$.
  \newline  It does not admit; $\phi \ngtr0$

\item \textbf{(123457F)!!}: $\dim(\Der)=10$, $\rank=0$, Betti Numbers $(2, 3, 4, 4, 3, 2, 1)$ ($\cong_{\RR} \ggo_{0.1}$!!).
  \newline Pre-Einstein der. $\phi=0$,
  \newline  It does not admit; $\phi \ngtr0$

\item \textbf{(123457H)!!}: $\dim(\Der)=10$, $\rank=0$, Betti Numbers $(2, 3, 4, 4, 3, 2, 1)$ ($\cong_{\RR} \ggo_{0.2}$!!).
  \newline Pre-Einstein der. $\phi=0$,
  \newline It does not admit; $\phi \ngtr0$.

\item \textbf{(123457H1)!!}: $\dim(\Der)=10$, $\rank=0$, Betti Numbers $(2, 3, 4, 4, 3, 2, 1)$ ($\cong_{\RR} \ggo_{0.2}$!!).
  \newline Pre-Einstein der. $\phi=0$,
  \newline It does not admit; $\phi \ngtr0$

\item\label{123457I} \textbf{(123457I) $[t\in\RR\setminus\{0,1\}]$}: $\dim(\Der)=10$, $\rank=1$, Betti Numbers $(2, 3, 4, 4, 3, 2, 1)$ ($\cong_{\RR} \ggo_{1.1(i_{t})}$)
  \newline Einstein der. $\phi=\frac{1}{7}\diag(1,2,3,4,5,6,7)$, $||\beta||^2=\frac{5}{7}\thickapprox0.7143$
  \newline It family admits nilsoliton metrics.

\item[(\ref{123457I}.1)] \textbf{(123457I) $[t=0]$}: $\dim(\Der)=10$, $\rank=1$, Betti Numbers $(2, 3, 4, 4, 3, 2, 1)$ ($\cong_{\RR} \ggo_{1.1(i_{0})}$)
  \newline Pre-Einstein der. $\phi=\frac{1}{5}\diag(1,2,3,4,5,6,7)$,
  \newline It does not admit.

\item[(\ref{123457I}.2)] \textbf{(123457I) $[t=1]$}: $\dim(\Der)=11$, $\rank=1$, Betti Numbers $(2, 3, 4, 4, 3, 2, 1)$ ($\cong_{\RR} \ggo_{1.1(i_{1})}$)
  \newline Pre-Einstein der. $\phi=\frac{1}{5}\diag(1,2,3,4,5,6,7)$,
  \newline It does not admit.
\end{enumerate}

\begin{example}\label{ex3}
For a final example, we consider the one-parameter family $(1357S)[t\in\RR\setminus\{0,1\}]$ given by $(\RR^7,\mu_{t})$ with
$$\mu_{t}=\left\{\begin{array}{l}
[e_1, e_2] = e_3, [e_1, e_3] = e_5, [e_1, e_5] = e_7,[e_1, e_6] = e_7, [e_2, e_3] = e_6, [e_2, e_4] = e_6, \\
{[e_2, e_5]} = e_7, [e_2, e_6] = te_7, [e_3, e_4] = e_7
\end{array}\right. .
$$
It is easily seen that $(1357S)[t]$ (with $t\neq0,1$) is a real form of the (complex) Lie algebra $(1357QRS1)[\lambda]$ with $\lambda:=\frac{t+\sqrt{t}}{t-\sqrt{t}}$. In the same manner as in the Example \ref{ex1}, we can see that the $\mathrm{GL}_7(\CC)$-orbit
of $(1357QRS1)[\lambda]$ is distinguished for the natural action of $\mathrm{GL}_7(\CC)$ on ${\Lambda^2(\CC^7)^*\otimes\CC^7}$
(by using Nikolayevsky's nice basis criterium in the complex case). Consequently, the $\mathrm{GL}_7(\RR)$-orbit
of $(1357S)[t]$ (with $t\neq0,1$) is distinguished; i.e. such family admits nilsoliton metrics (Theorem \ref{complex1}).

We fix the value of $t$, say $t=-3$. So $(1357S)[t=-3]$ is a real form of $(1357QRS1)[\frac{1}{2}-\frac{\sqrt {3}}{2}\sqrt{-1}] :=(\CC^7,\mu)$ with
$$
\mu=\left\{\begin{array}{l}
[{e_1},{e_2}]={e_3},[{e_1},{e_3}]={e_5},[{e_1},{e_4}]={e_6},[{e_1},{e_5}]={e_7},[{e_2},{e_3}]=-{e_6},[{e_2},{e_4}]={e_5},\\
{[{e_2},{e_6}]}= \left( \frac{1}{2}-\frac{\sqrt {3}}{2}\sqrt{-1}\right) {e_7},[{e_3},{e_4}]= \left( \frac{1}{2}+\frac{\sqrt {3}}{2}\sqrt{-1} \right) {e_7}
\end{array}\right. .
$$
To find a nilsoliton metric for $(1357S)[t=-3]$, we can find a distinguished point $\widetilde{\mu}$ in the $\mathrm{GL}_7(\CC)$-orbit of $\mu$ and then to study the real forms in the $\U(7)$-orbit of $\widetilde{\mu}$.

The point $\widetilde{\mu}$ can be found easily in much the same way as it was done for $(1357QRS1)[t=-1]$ in the Example \ref{ex1}. So we find $\widetilde{\mu}$ given by
$$
\widetilde{\mu}=\left\{\begin{array}{l} [{e_1},{e_2}]=\frac{\sqrt{13}}{13}{e_3},[{e_1},{e_3}]=\frac{\sqrt{13}}{13}{e_5},[{e_1},{e_4}]=\frac{\sqrt{26}}{26}{e_6},[{e_1},{e_5}]=\frac{\sqrt{13}}{13}{e_7},
{[{e_2},{e_3}]}=-\frac{\sqrt{13}}{13}{e_6},\\{[{e_2},{e_4}]}=\frac{\sqrt{26}}{26}{e_5},
[{e_2},{e_6}]=\frac{\sqrt {13}}{13} \left( \frac{1}{2}-\frac{\sqrt {3}}{2}\sqrt{-1} \right){e_7},
[{e_3},{e_4}]={\frac {\sqrt {26}}{26}} \left( \frac{1}{2}+\frac{\sqrt {3}}{2}\sqrt{-1} \right) {e_7}
\end{array}\right. .
$$
Since $(1357S)[t=-3]$ and $(\CC^{7},\widetilde{\mu})$ have the \comillas{same} pre-Einstein derivation, $\frac{4}{11}\diag(1,1,2,2,3,3,4)$, we can try to find
a distinguished point in the $\mathrm{GL}_7(\RR)$-orbit of $(1357S)[t=-3]$ by considering the real forms in the $\mathrm{G}$-orbit of $\widetilde{\mu}$, where
$$
\mathrm{G}=\left\{g \in \U(7) : g=\diag \left( \left(\begin{array}{ll} z_{1,1}&z_{1,2}\\ z_{2,1}&z_{2,2}\end{array}\right),
\left(\begin{array}{ll} z_{3,3}&z_{3,4}\\ z_{4,3}&z_{4,4}\end{array}\right),
\left(\begin{array}{ll} z_{5,5}&z_{5,6}\\ z_{6,5}&z_{6,6}\end{array}\right),
z_{7,7} \right)\right\}.
$$
Therefore, we find
$$
\widehat{\mu}=\left\{\begin{array}{l}
[{e_1},{e_2}]=\frac{\sqrt{13}}{13}{e_3},[{e_1},{e_3}]=\frac{\sqrt{13}}{13}{e_5},[{e_1},{e_4}]=-\frac{\sqrt{26}}{26}{e_5},[{e_1},{e_5}]=\frac{\sqrt{39}}{26}{e_7},[{e_1},{e_6}]=\frac{\sqrt{13}}{26}{e_7},\\
{[{e_2},{e_3}]}=\frac{\sqrt{13}}{13}{e_6},[{e_2},{e_4}]=\frac{\sqrt{26}}{26}{e_6},[{e_2},{e_5}]=\frac{\sqrt{13}}{26}{e_7},[{e_2},{e_6}]=-\frac{\sqrt{39}}{26}{e_7},[{e_3},{e_4}]=\frac{\sqrt{26}}{26}{e_7}
\end{array}\right. .
$$
where $(1357S)[t=-3]$ is isomorphic to $(\RR^7,\widehat{\mu})$ (over $\RR$) and the canonical inner product of $\RR^7$ defines a nilsoliton
metric on $(\RR^7,\widehat{\mu})$.

By a similar argument, we can to prove that (147E1)$[t>1]$, which is a real form of the (complex) Lie algebra (147E)$[\lambda]$ with
$\lambda=\left(\frac{( 1-\sqrt{t^2-1}\sqrt{-1})}{t}\right)^2$, admits nilsoliton metrics.
\end{example}

\begin{remark}
Recently, it has been mentioned in \cite[Proposition 3.2]{FINO1} that $\ngo_{11}=(1357S)[t=-3]$ does not admit any nilsoliton metric, as also $\ngo_{12}=(147E1)[t=2]$ and
$\ngo_{10}\cong_{\RR} 1.3(i)[t=1]$. In the above example we give a nilsoliton metric for $\ngo_{11}$, and by a similar argument we find
a nilsoliton metric for $\ngo_{12}$ given by the change of basis
$$
g:=\diag \left(
\left(\begin{array}{ll} 0 & \frac{\sqrt {3}}{2} \\ 1 & -\frac{1}{2}\end{array}\right),
1,-\frac{1}{4},
\left(\begin{array}{ll} \frac{\sqrt {3}}{12} & \frac{\sqrt {3}}{6}\\ -\frac{1}{4}&0\end{array}\right),
-\frac{\sqrt {3}}{48} \right)
$$
which defines
$$
\widetilde{\mu}=\left\{\begin{array}{l}
[{e_1},{e_2}]=\frac{\sqrt{3}}{6}{e_4},[{e_1},{e_3}]=-\frac{\sqrt{3}}{6}{e_6},[{e_1},{e_5}]=-\frac{\sqrt{3}}{12}{e_7},[{e_1},{e_6}]=\frac{1}{4}\,{e_7},\\
{[{e_2},{e_3}]}=-\frac{\sqrt{3}}{6}{e_5},[{e_2},{e_5}]=\frac{1}{4}\,{e_7},[{e_2},{e_6}]=\frac{\sqrt{3}}{12}{e_7},[{e_3},{e_4}]=-\frac{\sqrt{3}}{6}{e_7}
\end{array}\right. .
$$
where the canonical inner product of $\RR^7$ is a nilsoliton metric of $(\RR^{7},\widetilde{\mu})\cong_{\RR} \ngo_{12}$. The nilpotent
Lie algebra $1.3(i)[t=1]$ also admits a nilsoliton metric and it was proved in \cite[Example 2]{FERNANDEZ-CULMA1}. The nilpotent Lie algebra
$\ngo_{9}$, in fact, does not admit a nilsoliton metric and this can be proved by using Nikolayevsky's nice basis criterium ($\ngo_{9} = 1.1(iv)$; see \cite[Table 1.]{FERNANDEZ-CULMA1}).
\end{remark}


\end{document}